\documentclass[11pt]{article}
\usepackage{amsmath,amsfonts}
\usepackage{graphicx}
\usepackage{epsfig}

\newcommand{\R}{\mathbb{R}}
\newcommand{\Q}{\mathbb{Q}}
\newtheorem{theorem}{Theorem}
\newtheorem{lemma}{Lemma}
\newtheorem{remark}{Remark}
\newcommand{\N}{\mathbb{N}}
\newcommand{\M}{\mathbb{M}}
\newcommand{\Fc}{\mathcal{F}}
\newcommand{\Hc}{\mathcal{H}}
\newcommand{\E}{\mathsf{E}}
\newcommand{\ONE}{{\bf 1}}
\newcommand{\Pp}{\mathsf{P}}
\newcommand{\Tb}{\mathbb{T}}

\newcommand{\kap}{\kappa}
\newcommand{\eps}{\varepsilon}
\newcommand{\bpf}[1][Proof]{{\noindent {\sc #1: }}}
\newcommand{\epf}{{{\hspace{4 ex} $\Box$ \smallskip}}}
\newcommand{\tst}{\textstyle}

\author{Yuri Bakhtin\thanks{School of Mathematics, Georgia Tech, Atlanta GA, 30332-0160;  email:bakhtin@math.gatech.edu, 404-894-9235 (office phone), 404-894-4409(fax)}}
\date{}
\title{SPDE Approximation for Random Trees}

\begin{document}

\maketitle

\begin{abstract} We consider the genealogy tree for a critical branching process conditioned on non-extinction.
We enumerate vertices in each generation of the tree so that for each two generations one can define a monotone map describing the ancestor--descendant relation between their vertices. We show that under appropriate rescaling this family of monotone maps converges in distribution in a special
topology to a limiting flow of discontinuous monotone maps which can be seen as a continuum tree. This flow is a solution of an SPDE with respect to a Brownian sheet.

\end{abstract}

\section{Introduction}
In this paper we present a point of view at large random trees.
In \cite{Bakhtin-RSA} we considered Boltzmann--Gibbs distributions on rooted plane trees with bounded branching and proved that
as the order of the tree grows to infinity, these trees obey a certain thermodynamic limit theorem with a limit
given by an infinite Markov random tree that can be seen as the genealogy tree for a specially chosen critical branching processes conditioned on non-extinction, see~\cite{Aldous-Pitman:MR1641670},\cite{Kesten:MR871905},\cite{Kennedy:MR0386042}.
Appropriately rescaled generation sizes in the limiting tree admit an approximation
by a diffusion process. However, this view ignores all the interesting details concerning the complicated way different generations
of the random tree are connected to each other.

  The goal of this paper is to extend the diffusion approximation result and investigate the fine structure of the infinite random tree. 
The idea is to encode the tree via a stochastic flow of monotone maps.
We show that under the same rescaling, the flow of monotone maps associated to the tree converges in distribution to a 
monotone flow that can be viewed as a solution of an evolutionary SPDE with respect to a Brownian sheet. The possibility of such a result was hinted at in~\cite{Bakhtin-RSA}, and in this text we provide
a rigorous treatment of the SPDE, some interesting properties of its solutions, and the convergence theorem in an
appropriate metric space. The solution of SPDE that we construct can be called a continuum random tree, and
it appears to be a new type of continuum random trees, unknown in the literature, although some other stochastic flows
similar to ours have appeared in the literature (see, e.g., \cite{LeGall:MR1714707}) and there is a natural connection of our results to superprocesses 
that describe the evolution of mass generated by continual
branching particle systems conditioned on nonextinction, see~\cite{Evans-Perkins:MR1088825},\cite{Evans:MR1249698}. 
We describe this connection and the advantages of our approach in~Section~\ref{sec:superprocesses}. In short, the results of this paper
along with those of~\cite{Bakhtin-RSA} show what a typical large tree looks like if drawn on the plane.

{\bf Acknowledgment.} The author wishes to thank Tom Kurtz, Carl Mueller and Andrey Pilipenko for their comments on connection
of this work to the existing literature on superprocesses. He is grateful to NSF for partial support through CAREER grant DMS-0742424.

\section{Gibbs distributions on plane trees and thermodynamic limit}\label{sec:Gibbs}
In this section we describe the background in detail. Recall that plane (or, ordered) trees are rooted trees
such that subtrees at any vertex are linearly ordered.

We fix $D\in\N$ and introduce $\Tb_N=\Tb_N(D)$, the set of all plane trees on~$N$ vertices such that the branching number (i.e.\ the number of children) of each vertex  does not exceed $D$. 

We assume that every admissible branching number $i\in\{0,\ldots,D\}$ is assigned an energy $E_i\in\R$, and the energy of the tree $T\in\Tb_N$ is defined via
\[
 E(T)=\sum_{v\in V(T)} E_{\deg(v)}=\sum_{i=0}^D \chi_i(T)E_i,
\]
where $V(T)$ denotes the set of vertices of the tree $T$, $\deg(v)$ denotes the branching number of vertex~$v$, and $\chi_i(T)$ is the number of vertices of degree~$i$ in $T$. This energy function defines a local interaction between
the vertices of the tree since the energy contribution from a single vertex depends only on its nearest neighborhood
in the tree.

We fix an inverse temperature parameter $\beta>0$ and define a probability measure $\mu_N$ on $\Tb_N$ by
\[
 \mu_N\{T\}=\frac{e^{-\beta E(T)}}{Z_N},
\]
 where the normalizing factor (partition function) is
\[
Z_N=\sum_{T\in \Tb_N} e^{-\beta H(T)}.
\]

In~\cite{Bakhtin-RSA}, the limiting behavior of measures $\mu_N$ as $N\to\infty$ was studied. To recall the results 
of~\cite{Bakhtin-RSA}, we need more notation and terminology.

For each vertex $v$ of a tree $T\in\Tb_N$ its height $h(v)$ is defined as the distance to the root of $T$, i.e.\ the length of the shortest path connecting~$v$ to the root along the edges of $T$. The height of a finite tree is the maximum height of its vertices.

Each rooted plane tree can be uniquely encoded as a sequence of generations. By a generation we mean a 
monotone (nondecreasing) map $G:\{1,\ldots,k\}\to\N$, or, equivalently, the set of pairs $\{(i,G(i)):\ i=1,\ldots,k\}$
such that if $i_1\le i_2$ then $G(i_1)\le G(i_2)$.
We denote $|G|=k$ the number of vertices in the generation, and for any
$i=1,\ldots,k$, $G(i)$ denotes $i$'s parent number in the previous generation.  

For two generations~$G$ and~$G'$ we  
write $G\lhd G'$ and say that $G'$ is a continuation of $G$ if $G'(|G'|)\le |G|$. Each tree of height~$n$ can be
viewed as a sequence of generations 
\[1\lhd G_1\lhd G_2\lhd\ldots\lhd G_n\lhd 0,\] where $1\lhd G_1$ means 
$G_1(|G_1|)=1$ (the 0-th generation
consists of a unique vertex, the root) and $G_n\lhd 0$ means that the generation $n+1$ is empty. 
Infinite sequences $1\lhd G_1\lhd G_2\lhd\ldots$ naturally encode infinite trees.

For any plane tree $T$ and any $n\in\N$, $\pi_n T$ denotes the neighborhood of the root of radius $n$, i.e.\ the subtree of $T$ spanned by all vertices with height not exceeding $n$. 

For any $n$ and sufficiently large $N$, the map $\pi_n$ pushes the measure~$\mu_N$ on~$\Tb_N$ forward to the measure $\mu_N\pi_n^{-1}$ on $S_n$, the set of all trees with height~$n$.

We introduce $\rho$ and $C$ as a unique solution of
\begin{equation*}
\frac{1}{C}=\sum_{i=0}^De^{-\beta E_i}\rho^i=\sum_{i=0}^Die^{-\beta E_i}\rho^i,
\end{equation*}
and define
\begin{equation*}
 p_i=Ce^{-\beta E_i}\rho^i,\quad i=0,1,\ldots,D.
\end{equation*}
Then the vector $(p_i)_{i=0}^D$ defines a probability distribution on $\{0,1,\ldots,D\}$ with mean $1$. This
vector plays the role of a minimizer for the free energy (or large deviation rate function) associated to this model, see  \cite{Bakhtin-Heitsch-1:MR2415118} and~\cite{Bakhtin-Heitsch-2}.


Some of the results of 
~\cite{Bakhtin-RSA} are summarized in the next theorem. 

\begin{theorem} \label{th:therm_limit} \begin{enumerate}
\item There is a unique measure $P$ on infinite rooted plane trees with branching bounded by $D$ such that for any $n\in\N$
\[
 \mu_N \pi_n^{-1} \stackrel{TV}{\to} P\pi_n^{-1}.
\]
\item Measure $P$ defines a Markov chain on generations $(G_n)_{n\ge 0}$. 
The transition probabilities are given by
\begin{equation}
\label{eq:markov_transition}
\Pp\{G_{n+1}=g'|\ G_{n}=g\}=\begin{cases}\frac{|g'|}{|g|}p_{i_1}\ldots p_{i_{|g|}},& g\lhd g',\\
                             0,&\mbox{ otherwise,}
                            \end{cases}
\end{equation}
where $i_k$, $k=1,\ldots,|g|$ denotes the number of vertices in generation $g'$ that are children of $k$-th vertex in generation $g$.
\end{enumerate}
\end{theorem}
\begin{remark} Our notation differs from the notation used in \cite{Bakhtin-RSA}. It is an easy exercise to check that the~
r.h.s.\ of \eqref{eq:markov_transition} equals the expression given in~\cite{Bakhtin-RSA}.
\end{remark}
\begin{remark} The infinite Markov tree from Theorem~\ref{th:therm_limit} can be obtained as a genealogy tree of a critical branching process conditioned
on nonextinction, see~\cite{Aldous-Pitman:MR1641670},\cite{Kesten:MR871905}, and \cite{Kennedy:MR0386042}.
\end{remark}

Let us also recall a functional limit theorem for $X_n=|G_n|$, the process of the Markov infinite tree's generation sizes.
 We
introduce moments of the distribution $p$:
\begin{equation}
\label{eq:definition_of_B_n}
 B_n=\sum_{i=0}^Di^np_i,\quad n\in\N,
\end{equation}
and its variance
\[
 \mu=B_2-B_1^2=B_2-1.
\]
Let us now fix a positive time $T$ and define
\[
 Z_n(t)=\frac{1}{\mu n}(X_{[nt]}+\{nt\}(X_{[nt+1]}-X_{[nt]})),\quad t\in[0,T],
\]
the rescaled process of linear interpolation between values of $X_k$ given at whole times $k=0,1,\ldots$. 

\begin{theorem} The distribution of $Z_n$ converges weakly in uniform topology on $C([0,T])$ to the distribution
of a diffusion process with generator
\begin{equation}
\label{eq:generator_for_x*}
 Lf(x)=(f'(x)+\frac{1}{2}xf''(x))\ONE_{\{x\ge 0\}},
\end{equation}
emitted from $0$ at time $0$.
\end{theorem}
Notice that the limiting process can be viewed as a weak solution of an SDE
\[
 dX(t)=dt+\sqrt{X(t)}dW(t)
\]
on $\R_+$ (zero is entrance and non-exit singular boundary point for $\R_+$).

This result captures the asymptotics of a rather rough characteristic of the tree, the generation size. It ignores all the interesting details of the way the generations are connected to each other. In this paper, we are interested in taking these details into account,
thus refining the diffusion approximation.

From now on we consider the random genealogy trees defined above via the probability vector $(p_i)_{i=0}^D$. Although $D$ will
always be assumed to be finite, the results also hold true if one replaces the boundedness of $D$ with certain moment restrictions
on the distribution. 

\section{Limit theorems for finite partitions}\label{sec:limit_for_finite_partitions}
Let us take a positive number $z$, and for each $n\in\N$ consider the Markov process on generations (described in the previous section) originating at time $0$ 
with a population of $[nz]$ vertices. Let us take $m\in\N$ and choose some partition $0=x_0<x_1<\ldots<x_{m-1}<x_m=z$ of $[0,z]$.
Let us denote by $U(n,k,j)$ the size of the progeny of first $[nx_k]$ vertices in the initial population after $j$ steps.
Notice that $V(n,k,j)=U(n,k,j)-U(n,{k-1},j)$ describes the evolution of the progeny of disjoint subpopulations of the original population.
Denoting the rescaled interpolation by
\[
V_{n,k}(t)=\frac{1}{\mu n} \bigl(V(n,k,[nt])+\{nt\}(V(n,k,[nt]+1)-V(n,k,[nt]))\bigr),\quad t\in[0,T],
\] 
and following the same lines as in \cite{Bakhtin-RSA}, we obtain the following result:
\begin{theorem}\label{th:limit-V} Vector-valued process $(V_{n,k})_{k=1}^m$ converges (as $n\to\infty$) in distribution
in the uniform topology to a diffusion process $(V_k)_{k=1}^m$  on~$\R_+^m$ with 
initial data $V_k(0)=v_k=x_k-x_{k-1}$ and generator $L$, for all $C^2$-functions with compact support given by
\begin{equation}
\label{eq:generator_V}
 Lf(y)=\sum_{k=1}^m \frac{y_k}{|y|_1}\ONE_{\{y_k\ge 0\}}\partial_k f(y)+\frac{1}{2}\sum_{k=1}^m y_k\ONE_{\{y_k\ge 0\}}\partial^2_{kk}f(y),\quad y\in \R_+^m,
\end{equation}
where $|y|_1=y_1+\ldots+y_m$.
\end{theorem}
 Equivalently, the limiting process can be represented as a coordinatewise nonnegative solution of a system of SDEs:
\begin{align}
dV_k(t)&=\frac{V_k(t)}{|V(t)|_1}dt+\sqrt{V_k(t)}dW_k(t),\label{eq:V-time-noise}\\
 V_k(0)&=v_k,\quad k=1,\ldots, m,\label{eq:V-initial}
\end{align}
where $(W_k)_{k=1}^m$ are independent Wiener processes.
 
It is clear that this system has a unique strong solution until the first time at which $V_k=0$
for some~$k$. It is also easy to see (essentially from the Feller classification of singular points for 1-dimensional SDEs) that for each $k=1,\ldots,m$, the coordinate hyperplane $\{x_k=0\}$ is an absorbing set. Therefore, as soon as the system reaches one of the coordinate planes, we have to solve
a system of $m-1$ equations from that point on. Iterating this process we obtain a unique strong solution in $\R_+^m$.

In the same way, a sequence $(U_{n,k})_{k=1}^m$ defined by 
\[
U_{n,k}(t)=\frac{1}{\mu n} \bigl(U(n,k,[nt])+\{nt\}(U(n,k,[nt]+1)-U(n,k,[nt]))\bigr),\quad t\in[0,T],
\]
satisfies a functional limit theorem.

\begin{theorem} The vector-valued process $(U_{n,k})_{k=1}^m$ converges in distribution in uniform topology to a limiting 
 process $(U_{n})_{k=1}^m$ that can be represented as
\[
U_k(t)=V_1(t)+\ldots+V_k(t),\quad k=1,\ldots,m, 
\]
where $V$ is the limiting process from Theorem~\ref{th:limit-V}.
\end{theorem}
 Equivalently, the limiting process can be represented as a (nondecreasing in~$k$) solution of a system of SDEs:
\begin{align}
dU_k(t)&=\frac{U_k(t)}{U_m(t)}dt+\sum_{j=1}^k\sqrt{U_k(t)-U_{k-1}(t)}dW_j(t),\label{eq:U-time-noise}\\
 U_k(0)&=x_k,\quad k=1,\ldots, m,\label{eq:U-initial}
\end{align}
where $(W_j)_{j=1}^m$ are independent Wiener processes.

Although the theorems above provide interesting information about the behavior of the progeny of finitely many subpopulations, we do not give detailed proofs of these results, since we are actually aiming at more advanced ones (nevertheless, see Section~\ref{sbsc:characterization} for a proof
of a similar statement).

It is important to notice that the system \eqref{eq:U-time-noise},\eqref{eq:U-initial} (or, equivalently, \eqref{eq:V-time-noise},\eqref{eq:V-initial}) not only determines the evolution of process $U$ (respectively $V$) in time, but also has
a rich ``spatial'' structure which will be explored in later sections.

\section{SPDE approximation: first steps}\label{sec:spde_approx_first_steps}

In this section we give an informal view at the possibility of constructing a limiting continuum random tree represented as a solution
to a certain SPDE.

 In the approach described in the last section, for any partition $(x_k)_{k=0}^m$ 
we were able to find a probability space and a diffusion process defined on it that serves as a weak limit for the 
evolution of subpopulations in a discrete tree. Note that for any $t\ge 0$, the pre-limit processes $V_n,U_n$ as well as their
limits $U,V$ constructed above give rise to a random monotone map $x_k\mapsto U_k(t)$ defined on $(x_k)_{k=0}^m$ and the aforementioned
probability space. 

A drawback of this approach is that we need a separate probability space for each partition whereas
it would be natural to have just one probability space and a random flow of monotone 
maps $x\mapsto U(x,t)$ (defined on that space) that serves
all partitions simultaneously. That is, we would like to require that for any $m$ and any partition~$(x_k)_{k=0}^m$ of $[0,z]$, $(U(x_k,t)-U(x_{k-1},t))_{k=1}^m$ is a Markov process with
generator $L$ defined in~\eqref{eq:generator_V}.

To that end we utilize the spatial structure that has been mentioned in the end of Section~\ref{sec:limit_for_finite_partitions} and introduce a stochastic equation with respect to a Brownian sheet $W$ defined on Borel subsets of $\R_+^2=\{(x,t): x,t\ge 0\}$
and viewed as an orthogonal martingale measure (see \cite[Chapter 2]{Walsh-SPDE:MR876085} for a definition of Brownian sheet and stochastic integration with respect to orthogonal martingale measures):
\begin{align*}
dU(x,t)&= \frac{U(x,t)}{U(z,t)}dt+W([0,U(x,t)]\times dt),\\
U(x,t_0)&=x,\quad x\in[0,z],
\end{align*}
or, in integral form,
\begin{equation}
U(x,t)=x+\int_0^t \frac{U(x,s)}{U(z,s)}ds+\int_0^t\int_{\R}\ONE_{[0,U(x,t)]}(y) W(dy\times dt),
\label{eq:spde-1-integral} 
\end{equation}

The choice of this SPDE is a natural one, since the drift term and the quadratic covariation structure of the martingale part coincide with 
those of~\eqref{eq:U-time-noise}.

Our discussion of this SPDE will concentrate around the following questions:
 a natural definition of a solution, its existence and uniqueness, regularity properties, and the relation to the finite-dimensional 
system~\eqref{eq:U-time-noise},\eqref{eq:U-initial}.

First, let us take any $z>0$ and
$x=z$ in \eqref{eq:spde-1-integral}. The equation rewrites then as
\begin{equation}
U(z,t)=z+t+\int_0^t \int_{\R}\ONE_{[0,U(z,t)]}(y)W(dy\times dt).
\label{eq:spde-x*-integral} 
\end{equation}

One can obtain a unique strong solution of this equation using the standard Picard iteration scheme and mimicking
the proof of existence and uniqueness for an ordinary SDE. The only potential difficulty one may encounter is the behavior
of the solution near 0, but the solution is a Markov process
with generator $L$ given in~\eqref{eq:generator_for_x*} so that the solution always stays positive with probability 1.

Next, we can take any $z>0$ and a finite set $I$ of admissible initial points, i.e.,
\[\{0,z\}\subset I\subset [0,z].\] 
Mimicking the construction and the proof of the uniqueness of a positive solution of \eqref{eq:V-time-noise},\eqref{eq:V-initial}
we see that \eqref{eq:spde-1-integral} has a unique strong nondecreasing in $x\in I$ solution
 $(U(x,t))_{x\in I, t\ge 0}$.
Thus obtained $U$ satisfies a finite system of stochastic equations w.r.t.\ the Brownian sheet with probability one. 
It is clear that if $I'\supset I$ is a broader finite set of admissible initial conditions, then, with probability 1, for any $t\ge 0$, the map
$x\mapsto U(x,t), x\in I'$ is a monotone extension of the map $x\mapsto U(x,t), x\in I$.

In principle, there might be a problem with defining this map with probability 1 for an uncountable set $x\in[0,z]$. However, we can
use the monotonicity to tackle this difficulty. Let us take a countable set $J$ of admissible initial 
conditions such that $J$ is dense in $[0,z]$. Then, taking a sequence of finite sets increasing to $J$, and iteratively extending
the map $x\mapsto U(x,t)$ on each set of the sequence as described above, we are able to define this map on $J$ with probability 1, for all
$t\ge 0$. Moreover, this map is a.s.-monotone for all $t$. Therefore, this map can be extended uniquely by monotonicity to all points of $[0,z]$
except for at most countable set of discontinuities. Notice that for two different dense sets $J,J'$ of admissible initial conditions, with probability 1, the maps $U$ defined on $J$ and $J'$ are monotone continuations of each other. Therefore the monotone extensions
of these two maps agree at all points of $[0,z]$ except for a countable set of upward jumps. 

It is noteworthy that discontinuities in the form of upward jumps, or shocks, do occur with probability 1 as we shall see in
Section~\ref{sec:properties_of_solutions}. To deal with them we are
going to ignore the value of $U(x,t)$ at the jump, and work with equivalence classes of monotone functions coinciding at every continuity point. 
We proceed to introduce these equivalence classes in the next section. 


\section{Monotone graphs and flows}\label{sec:monotone_graphs}
The goal of this section is to introduce a space that will serve as a natural existence and uniqueness class of solutions of 
SPDE~\eqref{eq:spde-1-integral}.


Consider all points $z\ge 0$ and nonnegative nondecreasing functions $f$ defined on $(-\infty,z]$ such that
$f(x)=0$ for all $x<0$. Each of these
functions has at most countably many discontinuities. We say that two such functions $f_1:(-\infty,z_1]\to\R^+$, $f_2:(-\infty,z_2]\to\R^+$ 
are equivalent if $z_1=z_2$, $f_1(z_1)=f_2(z_2)$, and for each continuity point $x$ of $f_1$, $f_1(x)=f_2(x)$. Although
the roles of $f_1$ and $f_2$ seem to be different in this definition, it is easily seen to define a true equivalence relation.
The set of all classes of equivalence will be denoted by~$\M$. 
We would like to endow $\M$ with  a metric structure, and to that end we develop a couple of points of view. 

Sometimes, it is convenient to identify each element of $\M$ with its unique right-continuous
representative. Sometimes, it is also convenient to work with graphs. The graph of a monotone function $f$ defined on $(-\infty,z]$
is the set $G_f=\{(x,f(x))$, $x\le z\}$. For each 
discontinuity point $x$ of $f$ one may consider the line segment $\bar f(x)$ connecting points $[x,f(x-)]$ and $[x,f(x+)]$.  
The continuous version of $G_f$ is the union of $G_f$ and all segments $\bar f(x)$. It is often convenient to identify an element of $\M$
with a continuous version of its graph restricted to $\R_+^2$, and we shall do so from now on calling the elements of $\M$ monotone graphs.
Yet another way to look at monotone graphs is to think of them as monotone multivalued maps so that the image of each point is either
a point $f(x)$ or a segment $\bar f(x)$.


The Hausdorff distance between $\Gamma_1\in\M$ and $\Gamma_2\in\M$ is defined via
\[
 \rho(\Gamma_1,\Gamma_2)=\max\left\{\sup_{z_1\in \Gamma_1}\inf_{z_2\in \Gamma_2}|z_1-z_2|,\sup_{z_2\in \Gamma_2}\inf_{z_1\in \Gamma_1}|z_1-z_2|\right\}.
\]

\begin{lemma}\label{lm:M-is-Polish} $(\M,\rho)$ is a Polish (complete and separable) metric space. 
\end{lemma}
\bpf The separability follows since one can easily approximate any monotone graph by broken lines with finitely many 
rational vertices.

The metric space of compact subsets of the plane is complete with respect to Hausdorff metric. Therefore, to establish the completeness of $(\M,\rho)$  we need to prove that a limit of a convergent sequence of monotone graphs is necessarily a monotone graph. 

The limit set is obviously a connected one, with at most one point on each of the lines $\{(x_0,x_1):\ \alpha_0 x_0+ \alpha_1x_1=s\}$ for any $\alpha_0,\alpha_1>0$ and any $s\in\R$. Therefore it is a curve $\{\gamma_0(s),\gamma_1(s)\}$ parametrized by $s=\gamma_0(s)+\gamma_1(s)$,
and the monotonicity follows easily. \epf

We can derive another view at the space $\M$ from the proof above. We can look at each monotone graph in $\M$ in a rotated coordinate frame: then
we shall see that the set $\{(\gamma_0(s)+\gamma_1(s),\gamma_1(s)-\gamma_0(s))\}$ is a graph of 1-Lipschitz function.

We can introduce another distance based on this viewpoint. Take any $\Gamma_1,\Gamma_2\in\M$ and consider these curves in $\R_+^2$
in the rotated coordinate system as 1-Lipschits functions
$g_1$ and $g_2$
over $[0,s_1]$ and $[0,s_2]$ respectively. Define $s_*=s_1\wedge s_2$, and
\[
 \rho'(\Gamma_1,\Gamma_2)=|s_1-s_2|+\sup_{0\le s\le s_*} |g_1(s)-g_2(s)|.
\]

The proof of the following lemma an easy exercise.
\begin{lemma}\label{lm:equivalence_metrics} Distances $\rho$ and $\rho'$ are equivalent on $\M$.
\end{lemma}

The following is a useful criterion of convergence in $(\M,\rho)$.
\begin{lemma}\label{lm:convergence_criterion} A sequence of monotone graphs $(\Gamma_n)_{n\in\N}$ with
right-continuous representatives $f_n:[0,z_n]\to\R_+, n\in\N$, 
converges in $\rho$ to $\Gamma$ with right-continuous representative $f:[0,z]\to\R_+$ iff
\begin{enumerate}
 \item $z_n\to z$,\quad as $n\to\infty$;
 \item $f_n(z_n)\to f(z)$,\quad as $n\to\infty$;
 \item At each continuity point $x<z$ of $f$, $f_n(x)\to f(x)$, as $n\to\infty$.
\end{enumerate}
\end{lemma}
\bpf Suppose conditions 1--3 are satisfied. If $z=0,$ then the convergence is obvious. 
If $z>0$, take any $\eps>0$ and find continuity points $x_1<x_2<\ldots<x_m$ so that
$x_1<\eps$, $x_2-x_1<\eps$,\dots, $z-x_m<\eps$. Then take $n_0\in\N$ such that for all $n>n_0$, 
\begin{enumerate}
 \item[(i)]  for all $k=1,\ldots,m$,\quad $x_k<z_n$;
 \item[(ii)] for all $k=1,\ldots,m$,\quad $|f_n(x_k)-f(x_k)|<\eps$;
 \item[(iii)] $|z-z_n|<\eps$;
 \item[(iv)] $|f_n(z)-f(z_n)|<\eps$.
\end{enumerate}
It is easy to see that due to the monotonicity, for these values of $n$ we have $\rho(\Gamma_n,\Gamma)<2\eps$, and thus the convergence holds. 

Suppose now that the convergence holds. Conditions~1 and~2 of the lemma follow immediately by monotonicity. If condition~3 is violated, then
there is $\delta>0$ and a subsequence $(n')$ such that $|f_{n'}(x)-f(x)|>\delta$. We can choose $\eps>0$ so that $|y-x|<\eps$ implies
$|f(x)-f(y)|<\delta/2$. Therefore, for each $n'$, the point $(x,f_{n'}(x))$ is at least at distance $\eps\wedge (\delta/2)$ from any point
on $\Gamma$ which contradicts the assumption and completes the proof.
\epf

 The criterion of convergence provided by the above lemma allows to conclude immediately that any bounded set of monotone graphs is precompact.
Let us introduce $z_j(\Gamma)=\sup\{x_j:\ (x_0,x_1)\in\Gamma\}$, $j=0,1$.
\begin{lemma} A set $H\in\M$ is precompact in $\M$ iff
the set $\{(z_0(\Gamma),z_1(\Gamma)):\ \Gamma\in H\}$ is bounded in $\R^2$. 
\end{lemma}

This lemma can also be derived from Lemma~\ref{lm:equivalence_metrics} since the Lipschitz constant
for monotone graphs in the rotated coordinate system is bounded by~1.
 
Using Lemma~\ref{lm:convergence_criterion}, we can also prove the following statement, which we will not use but give here
for completeness:
\begin{lemma} For any $z>0$, the set $\M(z)$ of all monotone graphs $\Gamma\in\M$ with
$z_0(\Gamma)=z$ is a closed set. Convergence of monotone graphs in $\M(z)$ in distance $\rho$
is equivalent to essential convergence of their monotone right-continuous representatives which is
in turn equivalent to their convergence in Skorokhod topology on $D([0,z])$.
\end{lemma}

For two monotone graphs $\Gamma_1$ and $\Gamma_2$ with $z_1(\Gamma_1)=z_0(\Gamma_2)$,
we define their composition $\Gamma_2\circ\Gamma_1$ as the set of all pairs $(x_0,x_1)$ such that $(x_0,x_2)\in\Gamma_1$
and $(x_2,x_1)\in\Gamma_2$ for some $x_2$.

\bigskip

Let $T>0$ and $\Delta_{T}=\{(t_0, t_1): 0\le t_0\le t_1\le T\}$. We say that 
$(\Gamma^{t_0,t_1})_{(t_0,t_1)\in\Delta_T}$ is a (continuous) {\it monotone flow} on $[0,T]$ if the following properties
are satisfied:
\begin{enumerate}
 \item For each $(t_0,t_1)\in\Delta_T$, $\Gamma^{t_0,t_1}$ is a monotone graph.
 \item \label{pr:continuity} The monotone graph $\Gamma^{t_0,t_1}$ depends on $(t_0,t_1)$ continuously in $\rho$.
 \item \label{pr:z(t)}
For each $t\in[0,T]$, $\Gamma^{t,t}$ is the identity map on $[0,Z(t)]$ for some $Z(t)$. The function $Z$ is called the {\it profile}
of $\Gamma$.
 \item For any $(t_0,t_1)\in\Delta_T$, $z_0(\Gamma^{t_0,t_1})=Z(t_0)$, $z_1(\Gamma^{t_0,t_1})=Z(t_1)$, where
$Z$ is the profile of $\Gamma$.
\item\label{pr:agreement} If $(t_0,t_1)\in\Delta_T$ and $(t_1,t_2)\in\Delta_T$, then
$\Gamma^{t_0,t_2}=\Gamma^{t_1,t_2}\circ \Gamma^{t_0,t_1}$.  
\end{enumerate}

It is easy to check that the space $\M[0,T]$ of all monotone flows on $[0,T]$ is a Polish space if
equipped with uniform Hausdorff distance:
\begin{equation}
\label{eq:uniform_rho}
 \rho_{T}(\Gamma_1,\Gamma_2)=\sup_{(t_0,t_1)\in\Delta_T}\rho(\Gamma_1^{t_0,t_1},\Gamma_2^{t_0,t_1}).
\end{equation}

Property~\ref{pr:agreement} (consistency) implies that Property~\ref{pr:continuity} (continuity) has to be checked
only for $t_0=t_1$.

Our next goal is to introduce trajectories of individual points in the monotone flow.
Suppose $\Gamma\in\Delta_T$, and $Z$ is the profile of $\Gamma$. Let $U:\R^+\times\Delta_T\to\R_+$ satisfy the following properties:
\begin{enumerate}
 \item For any $(t_0,t_1)\in\Delta_T$,  the function $U(x,t_0,t_1)$ is monotone in $x\in[0,Z(t_0)]$. 
 \item For any $(t_0,t_1)\in\Delta_T$,  if $x\in[0,Z(t_0)]$, then  $(t_1,U(x,t_0,t_1))\in\Gamma^{t_0,t_1}$
 \item For all $x,t_0$, $U(x,t_0,t_1)$ is continuous in $t_1$.
\end{enumerate}

Then $U$ and $\Gamma$ are said to be compatible with
each other, and $U$ is said to be a {\it trajectory representation} of $\Gamma$.
Clearly, the monotonicity implies that, given $Z$, there is at most one monotone flow on $[0,T]$ compatible
with $U$. Moreover, it is sufficient to know a trajectory representation $U(x,t_0,t_1)$ for a dense set  of points
$x,t_0,t_1$  (e.g., rational points) to reconstruct the flow. 

Although a trajectory representation for a monotone flow  $\Gamma$ with profile $Z$ is not unique, there is a special
representation $U(x,t_0,t_1)$ that is right-continuous in $x\in[0,Z(t_0)]$ for every $t_1\ge t_0$:
\begin{equation}
 U(x,t_0,t_1)=\begin{cases}\sup\{y: (x,y)\in \Gamma^{t_0,t_1}\},&x\in[0,Z(t_0)]\\
              x,&x>Z(t_0). 
              \end{cases}
\label{eq:right_continuous_rep}
\end{equation}

The concrete way of defining $U(x,t_0,t_1)$ for $x>Z(t_0)$ is inessential for our purposes, and we often will simply
ignore points $(x,t_0,t_1)$ with $x>Z(t_0)$.

It is often convenient to understand a monotone flow as a triple $\Gamma=(\Gamma,Z,U)$, where $Z$ is the profile of $\Gamma$,
and $U$ is one of the trajectory representations of~$\Gamma$.

We are now ready to define a solution of our main equation~\eqref{eq:spde-1-integral}.



\section{Solution of the SPDE}\label{sec:spde_solution}

In this section, we define a solution of~\eqref{eq:spde-1-integral}, prove its existence, uniqueness, and the Markov
property.

Let $(\Omega,\Fc,\Pp)$ be a complete probability space with a filtration $(\Fc_t)_{t\ge0}$ satisfying the usual conditions 
(right-continuity and completeness), and let~$W$ be a Brownian sheet on $\R^2_+$ w.r.t.\ $(\Fc_t)_{t\ge 0}$, i.e.,
$W$ is a centered Gaussian random field indexed by Borel subsets of $\R^2_+$, such that
for any Borel set $A\subset\R_+$ with finite Lebesgue measure $|A|$, $(W([0,t]\times A))_{t\ge0},$ is
an $(\Fc_t)_{t\ge0}$-martingale with  $\langle W([0,\cdot]\times A)\rangle_t=t|A|$. We refer to~\cite{Walsh-SPDE:MR876085} for
more background on martingale measures.

We define a solution to~\eqref{eq:spde-1-integral} on $[0,T]$ 
as a random monotone flow $\Gamma:\Omega\to \M[0,T]$ such that with probability 1,  
the triplet $(\Gamma,Z,U)$ has the following properties for some trajectory representation $U$ of $\Gamma$: 
\begin{enumerate}
\item\label{it:measurability2} For any $t_0\ge0$ and any $t_1\ge t_0$, $\Gamma^{t_0,t_1}$ is measurable w.r.t.\ $\Fc_{t_1}$.
\item\label{it:measurability1} For any $x\ge 0$ and $(t_0, t_1)\in\Delta$, $U(x,t_0,t_1)$ is measurable w.r.t.~$\Fc_{t_1}$. 
\item \label{it:spde} For any $t_0\ge 0$, any $x\ge 0$, almost every $\omega\in\{x\le Z(t_0)\}$,
and all  $t_1\ge t_0$,
\begin{equation}
\label{eq:spde-with-X*}
 U(x,t_0,t_1)=x+\int_{t_0}^{t_1} \frac{U(x,t_0,t)}{Z(t)}dt+\int_{t_0}^{t_1}\int_{\R}\ONE_{[0,U(x,t_0,t)]}(y) W(dy\times dt).
\end{equation}
\end{enumerate}
Let us recall that the idea of this definition is to define a natural existence and uniqueness class of solutions by means of
consistent co-evolution of monotone graphs and trajectories of individual points.

\begin{theorem}\label{th:construction_of_solution} For any filtered probability space $(\Omega,\Fc,(\Fc_t)_{t\ge 0},\Pp)$ satisfying the 
usual conditions and a Brownian sheet $W$ w.r.t.\ $(\Fc_t)$, there is a unique solution $(\Gamma^{t_0,t_1})_{(t_0,t_1)\in\Delta}$ of 
equation~\eqref{eq:spde-1-integral} in the above sense. 
\end{theorem}
\bpf Let us construct a solution first. We begin with $Z$. Let us start with an auxiliary equation
\[
 Z(x,t)=x+t+\int_0^t\int_{\R}\ONE_{[0,Z(x,s)]}(y)W(dy\times ds)
\]
for $x>0$.
Mimicking the SDE methodology, it is easy to construct a unique strong solution of this equation defined up to a stopping time \[\tau_n=\inf\left\{t\ge0:\ Z(x,t)\le \frac{1}{n}\right\}.\] Letting $n\to\infty$ and noticing that $0$
is a no-exit singular point for the diffusion $Z(x,\cdot)$, we can extend this solution to the whole time semi-axis $\R_+$.
Using the uniqueness, we see that with probability 1, for all $m\in\N$ and $t\ge 0$, \[Z(1/m,t)\ge Z(1/(m+1),t),\]
so that the limit $Z(t)=\lim_{m\to\infty} Z(1/m,t)$ is well-defined. It is easy to check now that $Z(t)$ is a strong solution
of equation
\[
Z(t)=t+\int_0^t\int_{\R}\ONE_{[0,Z(s)]}(y)W(dy\times dt).
\]
It is also easy to check that this strong solution is a unique one by mimicking the proof of the uniqueness theorem of Yamada and 
Watanabe for SDEs, see Proposition~2.13 in \cite[Chapter 5]{Karatzas--Shreve})

For any two nonnegative rational numbers $x,t$, we define a stochastic process $(U(x,t_0,t_1))_{t_1\ge t_0}$
as follows. If $x>Z(t_0)$, we set $U(x,t_0,t_1)=x$; if $x\le Z(t_0)$, we set $U(x,t_0,\cdot)$ to be
the unique strong solution of equation~\eqref{eq:spde-with-X*} on $[t_0,\infty)$ (the uniqueness can be established easily
using methods borrowed from the theory of SDEs). 

For each $t_1\ge t_0$ the resulting map $x\mapsto U(x,t_0,t_1)$ is nondecreasing on $[0,Z(t_0)]$: 
if $x_1\le x_2$ then $U(x_1,t_0,t_1)\le U(x_2,t_0,t_1)$. In fact, if the opposite inequality holds for some $t_1>t_0$ 
then $U(x_1,t_0,t')= U(x_2,t_0,t')$
for some $t'\in[t_0,t_1]$, and the uniqueness implies that $U(x_1,t_0,\cdot)$ and $U(x_2,t_0,\cdot)$ coincide after $t'$.

Now, for any $(t_0,t_1)\in\Delta$ with rational $t_0$, the (random) map $x\to U(x,t_0,t_1)$ defined above uniquely determines
a monotone graph $\Gamma^{t_0,t_1}$. Using the continuity of trajectories $U(x,t_0,t_1)$ and their monotonicity property,
it is easy to show that $\Gamma^{t_0,t_1}$ depends continuously on rationals $t_0$ and~$t_1$. Therefore, we can extend $\Gamma^{t_0,t_1}$
by continuity to all values $(t_0,t_1)\in\Delta$. 

Knowing monotone functions $U(x,t_0,t_1)$ for rational values of $x$ and $t_0$, and all $t_1\ge t_0$ allows to define
\[
U(x,t_0,t_1)=\inf\{U(y,t,t_1):\ t,y\in\Q,\ t<t_0,\ y\in[0,Z(t)],\ U(y,t,t_0)>x\},
\]
for all other values of $t_0$ and $x\in[0,Z(t_0)]$, and it is easy to verify that thus defined random 
field~$U$ satisfies the requirement of the Theorem.

To prove the uniqueness, it suffices to notice that the result of each step 
in the above construction of is a.s.-uniquely defined.
\epf
 
The proof of the theorem provides a construction of a unique solution of SPDE~\eqref{eq:spde-1-integral}.
In the next section, we study some properties of the solution. 
\begin{remark} The set of all pairs of rationals $x,t$ 
in the proof of Theorem~\ref{th:construction_of_solution}
can be replaced by any countable and dense set in $\R_+^2$ 
\end{remark}

The following statement follows directly from the Markov property of $U(x,t_0,t_1)$ in $t_1$:
\begin{lemma}[Markov property] For any $t_0\ge 0$, $(\Gamma^{t_0,t_1}, t_1\ge t_0)$ is a Markov process.
\end{lemma}

\section{Discontinuities of the solution}\label{sec:properties_of_solutions}

The reason for introducing monotone flows the way we did in Sections~\ref{sec:monotone_graphs} 
and~\ref{sec:spde_solution} is the presence of discontinuities of solutions with respect to $x\in[0,Z(t)]$. If not
for these discontinuities our analysis would have been easier, and in this section we 
will show that they are, in fact, intrinsic to the solution, thus justifying 
our choice to work with monotone graphs.

Let us call a straight line segment connecting two points $(x_1,y_1)$ and~$(x_2,y_2)$ on the plain vertical if $x_1=x_2$.

\begin{lemma}\label{lm:shocks}
Let $0\le t_0<t_1$. Then, with probability 1, $\Gamma^{t_0,t_1}$ contains nondegenerate vertical segments,
i.e., any monotone function representing $\Gamma^{t_0,t_1}$ 
is discontinuous, with shocks associated to these vertical segments.
\end{lemma}
\bpf
To prove this lemma we need the following standard result (see e.g. the proof of Theorem~26 in \cite[Chapter 2]{Protter:MR2020294}
for this result in the context of realizations of stochastic processes).

\begin{lemma}
 If $f:[0,z]\to\R$ is a bounded variation function then its quadratic variation $Q(f)$ 
defined by
\[
Q(f)=\lim_{n\to\infty}\sum_{i=0}^{n-1} \left(f(z(i+1)/n)-f(zi/n)\right)^2 
\]
satisfies:
\[
 Q(f)=\sum_{x\in \Delta(f)} (f(x+)-f(x-))^2,
\]
where
and $\Delta(f)$ is the set of all discontinuity points of $f$. In particular, $f$ is continuous on $[0,z]$ iff $Q(f)=0$.
\end{lemma}

Let us compute the quadratic variation of $U(\cdot,t_0,t_1)$. For $n\in\N$ we introduce (omitting dependence on $n$)  \[x_k=\frac{k}{n}Z(t_0),\quad k=0,\ldots,n,\]
 and
\[V_k(t_1)=U(x_k,t_0,t_1)-U(x_{k-1},t_0,t_1),\quad k=1,\ldots,n.\]
It\^o's formula implies
\begin{align*}
 dV_k^2(t_1)=&2V_k(t_1)dV_k(t_1)+V_k(t_1) dt_1\\
&=\frac{2 V_k^2(t_1)}{Z(t_1)}dt_1+2V_k(t_1)W([U(x_{k-1},t_0,t_1),U(x_{k},t_0,t_1)]\times dt_1)\\
&+V_k(t_1)dt_1.
\end{align*}
Let \[Q_n(t_1)=V_1^2(t_1)+\ldots+V_n^2(t_1).\] Then $Q_n(t_0)=(Z(t_0))^2/n$, and 
\begin{align*}
 Q_n(t_1)=&\frac{(Z(t_0))^2}{n}+2\int_{t_0}^{t_1}\frac{Q_n(t)}{Z(t)}dt\\&+2\sum_{k=1}^n \int_{t_0}^{t_1}V_k(t)W([U(x_{k-1},t_0,t),U(x_{k},t_0,t)]\times dt)\\&+\int_{t_0}^{t_1} Z(t) dt.
\end{align*}
Let us define $Q(t)=\lim_{n\to\infty}Q_n(t)$ and $\nu=\inf\{t>t_0: Q(t)>0\}$. If $\nu>t_0$, then taking the limit as $n\to\infty$ in both sides of the equation above at $t_1=\nu$, we see that 
all terms converge to zero  except for 
$\int_{t_0}^\nu Z(t)dt$. To obtain the convergence to zero for the martingale stochastic integral term, it is sufficient to see that the quadratic variation in time
given by
\begin{align*}
\left\langle\sum_{k=1}^n \int_{t_0}^{\cdot}V_k(t)W([U(x_{k-1},t_0,t),U(x_{k},t_0,t)]\times dt)\right\rangle_\nu=
  \int_{t_0}^{\nu} \sum_{k=1}^n V_k^3(t)dt,
\end{align*}
converges to zero as $n\to\infty$.

Since $\int_{t_0}^\nu Z(t)dt$ is a strictly positive random variable that does not depend on $n$, we obtain a contradiction which shows that $\nu=t_0$, so that $Q(t)>0$ for any $t>t_0$. \epf

\section{Convergence}

This section is the central part of the paper. Here we 
prove that the discrete infinite random genealogy tree $\tau$ converges in distribution 
in an appropriate sense under appropriate rescaling to the continuum random tree given by
the monotone flow described in the previous section. 

\subsection{Trees as monotone flows. Main result.}
To start with, we introduce a monotone flow associated with a realization of the infinite discrete tree $\tau$.
This procedure is analogous to that of linear interpolation for discrete time random walks leading to
Donsker's invariance principle.
 
We begin with an imbedding of the tree in the plane. Recall that there are $X_n\ge 1$ vertices in the $n$-th generation of the tree.
 For $i\in\{1,\ldots,X_n\}$, the $i$-th vertex of $n$-th generation is represented by
the point $(n,i-1)$ on the plane. The parent-child relation between two vertices of the tree
is represented by a straight line segment connecting the representations of these vertices.

Besides these ``regular'' segments, we shall need some auxiliary segments that are not an intrinsic part of the tree but will be used in representing
the discrete tree as a continuous flow.  Suppose a vertex $i$ in $n$-th generation has no children. Let $j$ be the maximal vertex in generation 
$n+1$ among those having their parents preceding $i$ in generation $n$.
Then an auxiliary segment of type~I connects the points $(n,i-1)$ and $(n+1,j-1)$. If vertex $1$ in $n$-th generation has no children,
points $(n,0)$ and $(n+1,0)$ are connected by an auxiliary segment of type II. Auxiliary segments of type~III connect points $(n,i-1)$ and $(n,i)$ for $1\le i\le X_n-1$.

Every bounded connected component of the complement to the union of the above segments on the plane is
either a parallelogram with two vertical sides of length 1, or a triangle with one vertical side of length 1. One can treat
both shapes as trapezoids (with one of the parallel sides having zero length in the case of triangle).

\begin{figure}
\centering
\epsfig{file= 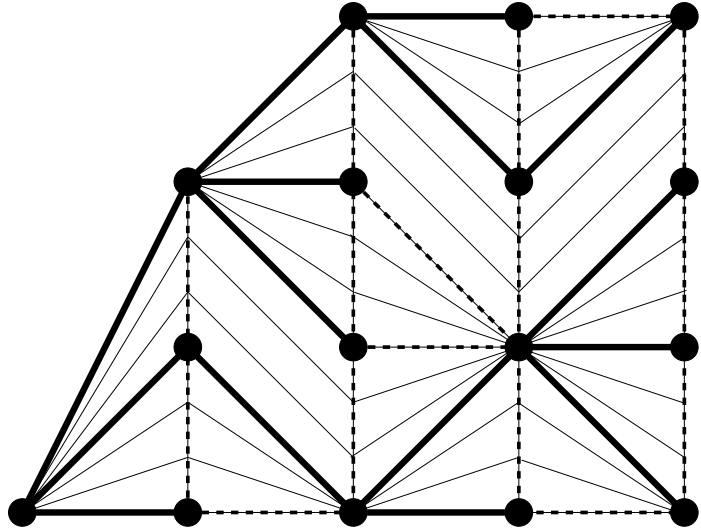,height=7cm}
\caption{Construction of the continuous monotone flow.}
\label{fig:construction}
\end{figure}

For each trapezoid, we shall establish a bijection with the unit square and define
the monotone flow to act along the images of the ``horizontal'' segments of the square. A graphic illustration of the construction 
is given on Figure~\ref{fig:construction}, and we proceed to describe it precisely.

Each trapezoid $L$ of this
family has vertices $g_{0,0}=(n,i_{0,0}),g_{0,1}=(n,i_{0,1})$, $g_{1,0}=(n+1,i_{1,0}), g_{1,1}=(n+1,i_{1,1})$, where $i_{0,1}-i_{0,0}\in\{0,1\}$ and $i_{1,1}-i_{1,0}\in\{0,1\}$.

Then, for every $\alpha\in(0,1)$ we define
\[
g^L_m(\alpha)=g_{m,0}+\alpha (g_{m,1}-g_{m,0}),\quad m=0,1,  
\]
and
\[
g^L(\alpha,s)=g^L_0(\alpha)+s (g^L_1(\alpha)-g^L_0(\alpha)),\quad s\in(0,1). 
\]
This definition introduces a coordinate system in $L$, i.e., a bijection between~$L$ and the unit square $(0,1)^2$.
We are going to use it to construct the monotone map associated with the tree for times $t_0,t_1$ assuming
that there is $n\in\{0\}\cup\N$ such that $n< t_0 \le t_1 <n+1$.
Let us take any $x$ such that $(x,t_0)$ belongs to one of the trapezoids $L$. Then there is a unique number
$\alpha(x,t_0)\in(0,1)$ such that $g^L(\alpha(x,t_0),\{t_0\})=x$, where $\{\cdot\}$ denotes the fractional part. We can define $g^{t_0,t_1}(x)=g^L(\alpha(x,t_0),\{t_1\})$.
This strictly increasing function can be consistently and uniquely extended by continuity to points $x$ such that $(x,t_0)$ belongs
either to a regular segment in the tree representation or an auxiliary segment of type~I or~II. This 
function $g^{t_0,t_1}$
also uniquely defines a monotone graph $\tilde \Gamma^{t_0,t_1}=\tilde \Gamma^{t_0,t_1}(\tau)$ depending continuously on $t_0,t_1$.
Next, if we allow $t_0$ and $t_1$ to take values $n$ and $n+1$, then we can construct the associated monotone graph
as the limit in $(\M,\rho)$ of the monotone graphs associated to the increasing functions defined above
(as $t_0\to n$ or $t_1\to n+1$). Notice that the resulting monotone graphs may have intervals of constancy
and shocks (i.e., contain horizontal and vertical segments). Now we can take any $(t_0,t_1)\in\Delta_\infty$ and define
\[
\tilde \Gamma^{t_0,t_1}=\tilde \Gamma^{[t_1],t_1}\circ\tilde \Gamma^{[t_1]-1,t_1}\circ\ldots\circ\tilde \Gamma^{[t_0]+1,[t_0]+2}\circ
\tilde \Gamma^{t_0,[t_0]+1}. 
\]
which results in a continuous monotone flow $(\tilde \Gamma^{t_0,t_1}(\tau))_{(t_0,t_1)\in\Delta_\infty}$.

To state our main result we need to introduce a rescaling of this family. For every $n\in\N$, we define
\begin{equation}
\label{eq:pre-limit-random-flows}
 \Gamma^{t_0,t_1}_n(\tau)=\left\{\left(\frac{x}{\mu n},\frac{y}{\mu n}\right): (x,y)\in\tilde \Gamma^{nt_0,nt_1}(\tau)\right\},\quad (t_0,t_1)\in\Delta_\infty.
\end{equation}
For each $T>0$ we can consider the uniform distance $\rho_{T}$ on monotone flows in $\M[0,T]$  introduced in~\eqref{eq:uniform_rho},
and define the locally uniform (LU) metric on~$\M[0,\infty)$ via
\[
 d(\Gamma_1,\Gamma_2)=\sum_{m=1}^\infty 2^{-m}(\rho_{m}(\Gamma_1,\Gamma_2)\wedge 1).
\]

\begin{theorem}\label{th:main_convergence} As $n\to\infty$, the random field $(\Gamma_n ^{t_0,t_1}(\tau))_{(t_0,t_1)\in\Delta}$
converges in distribution in LU metric to the stochastic flow of monotone graphs solving equation~\eqref{eq:spde-1-integral}.
\end{theorem}

The rest of the section is devoted to the proof of this theorem.
First, it is sufficient to prove the convergence in distribution in $\rho_{m}$ for all $m$. We take $m=1$ without loss of generality.
Due to the classical Prokhorov theorem it is sufficient to demonstrate two facts:
\begin{enumerate}
\item The sequence of distributions of $\Gamma_n(\tau)$ is tight w.r.t.\ $\rho_{1}$.
\item Any limit point for the sequence of distributions coincides with the distribution of the flow generated by SPDE.
\end{enumerate}

Subsection \ref{sbsc:tightness} is devoted to the proof of the first statement and subsection~\ref{sbsc:characterization} to the second one.

\subsection{Tightness}\label{sbsc:tightness}
We have to show that for any $\eps>0$ there is a set $K\subset\M[0,1]$  that is compact in $d_1$
and $\Pp\{\Gamma_n\in K\}>1-\eps$ for all $n$ (from now on we do not distinguish between $\Gamma_n$ and its restriction on $[0,1]$.)

The construction below depends on the values of constants $\alpha>0,\beta>\gamma>0$ and $b>0$. For $m\in\N$, consider a set 
\[
R_m=\left\{(k2^{-m},j2^{-m}):\ 0\le k \le 2^m;\ 0\le j\le 2^{m(1+\alpha)}\right\}, 
\]
and also define $R=\bigcup_m R_m$.
For a function $f$ we denote its $\beta$-H\"older constant by $\Hc_\beta(f)$:
\[
 \Hc_\beta(f)=\sup_{t_1,t_2:t_1\ne t_2}\frac{|f(t_1)-f(t_2)|}{|t_1-t_2|^\beta}.
\]
Let $K_m$ be the set of all monotone flows $(\Gamma,Z,U)$ in $\M[0,1]$ such that $\sup_{t}Z(t)\le 2^{m\alpha}$,
$\Hc_\beta(Z)\le b2^{m\gamma}$, and for each $(t_0,x)\in R_m$
satisfying $x<Z(t_0)$,
$\Hc_\beta(U(x,t_0,\cdot))\le b 2^{m\gamma}$. 
Denote $K=\bigcap_{m=m_0}^\infty K_m$. Note that we suppress the dependence of $K$ on $b$ and $m_0$.

\begin{lemma} For any $b$ and $m_0$, the set $K$ is precompact in $\M[0,1]$.\end{lemma}

\bpf  For any sequence $(\Gamma_n,Z_n,U_n)_{n\in\N}$ in $K$, we must show that it contains a convergent subsequence. 

Using the H\"older continuity, the Arcela--Ascoli compactness criterion, and the classical diagonal procedure,
we can extract a subsequence $n'$ such that for each $m\in\N$ and each
each $(t_0,x)\in R_m$, $U_{n'}(x,t_0,\cdot)$ converges uniformly to a limiting trajectory $U_\infty(x,t_0,\cdot)$ with
$\Hc_\beta(U_\infty(x,t_0,\cdot))\le b 2^{m\gamma}$, and $Z_{n'}$ converges uniformly to a limiting trajectory
$Z_{\infty}$.

We need to show that there is a unique continuous monotone flow $\Gamma_\infty$ compatible with these trajectories and that
$\Gamma_{n'}$ converges uniformly to $\Gamma_\infty$. Let us construct $\Gamma_\infty=(\Gamma_\infty,Z_\infty,U_\infty)$.
Notice that $Z_\infty$ is already at our disposal. 
Take any point
$(t_0,x)\in[0,1]\times \R_+$ such that $x<Z_\infty(t_0)$. For any $m$, we set  \[k=k(m)=[2^mt_0],\] 
\[
j^-=j^-(m)=[2^{m}(x-2^{m\gamma}\cdot (2^{-m})^\beta)],
\]
and
\[
j^+=j^+(m)=[2^{m}(x+2^{m\gamma}\cdot (2^{-m})^\beta)]+1.
\]
Then, 
\[j^+ 2^{-m}-j^-2^{-m}<2(2^{-m}+2^{m\gamma}\cdot 2^{-m\beta}), \]
and 
\[j^-2^{-m}+2^{m\gamma}\cdot (2^{-m})^\beta\le x \le j^+2^{-m}-2^{m\gamma}\cdot (2^{-m})^\beta.\]
For sufficiently large $m$,  $j^{+}2^{-m}< Z_\infty(k2^{-m})$, so that 
$U_\infty(j^\pm 2^{-m},k2^{-m},t_0)$ are already defined. Using the bound on the H\"older constant, we get
\begin{multline*}
x-2^{-m}-2\cdot 2^{m(\gamma-\beta)}\le U_\infty(j^-2^{-m},k2^{-m}, t_0)\le x\\ \le U_\infty(j^+2^{-m},k2^{-m}, t_0)\le x+2^{-m}+2\cdot2^{m(\gamma-\beta)}.
\end{multline*}
For any $t_1\ge t_0$, we define
\[
U^-(x, t_0,t_1)=\sup_m U_\infty(k(m)2^{-m},j^-(m)2^{-m}, t_1),
\]
and
\[
U^+(x,t_0,t_1)=\inf_m U_\infty(k(m)2^{-m},j^+(m)2^{-m}, t_1).
\]
It is clear that $U^-(x, t_0,t_1)\le U^+(x,t_0,t_1)$ and if $x<y$ then  $U^+(x, t_0,t_1)\le U^-(y,t_0,t_1)$.
In particular, both quantities are monotone in first argument and there is a monotone graph $\Gamma_\infty^{t_0,t_1}$ such that 
$U^-(x, t_0,t_1)$ and $U^+(x,t_0,t_1)$ are its left continuous and, respectively, right continuous representatives.
Also any monotone function compatible with $(U(x,t_0,t_1), (t_0,x)\in R)$, must represent
the same monotone graph.

The resulting family of monotone graphs $(\Gamma_\infty^{t_0,t_1})$ is easily seen to satisfy all the properties
in the definition of a monotone flow. Let us prove only Property 2, the continuity
of $\Gamma_\infty^{t_0,t_1}$ in $t_0$ and $t_1$.
It is sufficient to show that $\Gamma_\infty^{t_0,t_1}$ converges to the graph of the identity map on $[0,Z_\infty(t)]$ as
$t_0$ and~$t_1$ converge to $t$ from below and above respectively. 
To see the latter, we imbed a small time interval $(t_0,t_1)$ into some dyadic interval
$[k2^{-m},(k+1)2^{-m}]$ with a large  $m$. On $[k2^{-m},(k+1)2^{-m}]$, the monotone flow  $\Gamma_\infty$ displaces all dyadic points with denominator $2^{-m}$ by at most
$2^{m(\gamma-\beta)}$. Since $2^{m(\gamma-\beta)}\to 0$ as $m\to\infty$, the continuity follows.

Now we shall prove the convergence to this monotone flow $\Gamma_\infty$. 
For any $\eps>0$, we will show that for sufficiently
large values of $n'$, $d_1(\Gamma_\infty,\Gamma_{n'})<\eps$. First, for any $m\in\N$, we can find $n_0(m)$ such that for all $n'>n_0(m)$ and all $(t_0,x)\in R_m$,
$\|U_n(x,t_0,\cdot)-U_\infty(x,t_0,\cdot)\|<\eps/2$, where $\|\cdot\|$ denotes the sup-norm.
For any $(t_0,x)$ choose points $k,j^-,j^+$ as above and denote
\begin{align*}
 y^\pm_{n'}&=U_{n'}(j^\pm 2^{-m},k2^{-m},t_0),\\
y^\pm_{\infty}&=U_{\infty}(j^\pm 2^{-m},k2^{-m},t_0).\\
\end{align*}
Then
\begin{align*}
x-2^{-m}-2\cdot 2^{m(\gamma-\beta)}&\le y^-_{n'} \le x\le y^+_{n'}\le x+2^{-m}+2\cdot 2^{m(\gamma-\beta)},\\
 x-2^{-m}-2\cdot 2^{m(\gamma-\beta)}&\le y^-_\infty \le x\le y^+_\infty\le x+2^{-m}+2\cdot 2^{m(\gamma-\beta)}.
\end{align*}
The image of  $x$ under $\Gamma_{n'}^{t_0,t_1}$ (viewed as a multivalued map) for every $t_1>t_0$ is contained
in the image of $[y^-_{n'},y^+_{n'}]$ under $\Gamma_{n'}^{t_0,t_1}$. Therefore, the definition of $y^\pm_{n'}$ and
$y^\pm_{\infty}$ implies that
if $n'>n_0(m)$, then any point in the image of $[y^-_{n'},y^+_{n'}]$ under $\Gamma_{n'}^{t_0,t_1}$ is at most at $\eps/2$ from the image of $[y^-_{\infty},y^+_{\infty}]$ under~$\Gamma_\infty^{t_0,t_1}$. On the plane, the distance between any point $(x,x')\in\Gamma_{n'}^{t_0,t_1}$ and the graph
$\Gamma_\infty^{t_0,t_1}$ restricted to $[y^-_{\infty},y^+_{\infty}]$ does not exceed $2^{-m}-2\cdot 2^{m(\gamma-\beta)}+\eps/2$. So it is sufficient to choose
$m$ so that $2^{-m}-2\cdot 2^{m(\gamma-\beta)}<\eps/2$ to have that distance bounded by~$\eps$. Similarly, the same bound holds true
for the distance between any point $(x,x')\in\Gamma_\infty^{t_0,t_1}$ and the restriction of $\Gamma_{n'}^{t_0,t_1}$ to $[y^-_{\infty},y^+_{\infty}]$. We conclude that
for the above choice of $m$ and for any $n'>n_0(m)$, the uniform distance between $\Gamma_{n'}$ and $\Gamma_\infty$ is bounded by $\eps$, and the
convergence implying the precompactness of $K$ follows.
\epf

Our next goal is to show that for any $\eps>0$, the numbers $\alpha,\beta,\gamma,b$ can be chosen so that
$\Pp\{\Gamma_n\in K\}>1-\eps$ for all $n$.

\begin{lemma} \label{lm:holder_from_moment1}Suppose there are positive numbers $C,r$ such that for all $m\in\N$, all $(t_0,x)\in R_m$,
 a continuous process $U(x,t_0,\cdot)$ satisfies
\begin{equation}
\label{eq:moment_estimate}
 \E |U(x,t_0,t)-U(x,t_0,t')|^r \le C|t-t'|^{r/2},\quad t_0\le t\le t'.
\end{equation}
If
\begin{equation}\label{eq:r_should_be_large}
r/2-1>r\beta, 
\end{equation}
then there is a constant $C_1(C,\beta,r)$ such that for all $m\in\N$ and any $b>0$,
\[
 \Pp\left\{\sup_{(t_0,x)\in R_m} \Hc_\beta(U(x,t_0,\cdot)) > b2^{\gamma m}\right\}\le  C_1(C,\beta,r)\frac{2^{(2+\alpha-r\gamma)m}}{b^r}.
\]
\end{lemma}
\bpf Let us fix $x,t_0$ and (for brevity) denote $U(t)=U(x,t_0,t)$, 

Let us estimate $\Pp\{\Hc_\beta(U)\ge c\}$ for a $c\ge2$. 
If $\Hc_\beta(U)\ge c$, then there are times $t<t'$ such that
\[
 |U(t')-U(t)|>c(t'-t)^\beta.
\]
We can find $m$ such that $2^{-m}\le t'-t <2^{-m+1}$. Then 
\begin{equation}
\label{eq:increment_big}
 |U(t')-U(t)|> c 2^{-m\beta}.
\end{equation}
There is an integer $j$ such that 
$j2^{-m}\in[t,t']$. We can find sequences $(\kap_k)_{k=m}^\infty$ and $(\kap'_k)_{k=m}^\infty$ such that each $\kap_k$ and $\kap'_k$
is 0 or 1 and
\[
 t'-j{2^{-m}}=\sum_{k=m}^\infty\kap'_k2^{-k},
\]
\[
 j2^{-m}-t=\sum_{k=m}^\infty\kap_k2^{-k}.
\]
Assuming that on each respective dyadic interval of length $2^{-k}$, the increment of $U$ does not exceed
$(1-2^{-\beta})c2^{-\beta k-1}$ in absolute value, we have
\begin{align*}
 |U(t')-U(t)|&\le (1-2^{-\beta})c2^{-1}\left( \sum_{k=m}^\infty\kap'_k2^{-\beta k}+\sum_{k=m}^\infty\kap_k2^{-\beta k}\right)\\
&\le (1-2^{-\beta})c2^{-1}\cdot 2 \sum_{k=m}^\infty 2^{-\beta k}\\
&\le c 2^{-m\beta},
\end{align*}
which contradicts~\eqref{eq:increment_big}. Therefore, denoting $\nu(\beta)=(1-2^{-\beta})2^{-1}$, we see that
for some $k$ there is a dyadic interval $[l2^{-k},(l+1)2^{-k}]$ 
of length $2^{-k}$ such that 
\[
|U(l2^{-k})-U((l+1)2^{-k})|>\nu(\beta)c2^{-\beta k}.
\]
We use Markov's inequality and condition~\eqref{eq:moment_estimate} to see that an upper bound for the probability 
of this event is:
\begin{align*}
\sum_{k=0}^{\infty} 2^k \Pp\left\{ |U(l2^{-k})-U((l+1)2^{-k})|>\nu(\beta)c2^{-\beta k}\right\}\\
\le \frac{1}{c^{r}\nu(\beta)^r}\sum_{k=0}^{\infty} 2^k \E|U(l2^{-k})-U((l+1)2^{-k})|^r 2^{r\beta k}\\
\le  \frac{C}{c^{r}\nu(\beta)^r}\sum_{k=0}^{\infty} 2^k 2^{-kr/2} 2^{r\beta k}
\le \frac{C_1(C,\beta,r)}{c^r},
\end{align*}
where 
\[
C_1(C,\beta,r)=\frac{C}{\nu(\beta)^r}\sum_{k=0}^{\infty} 2^k 2^{-kr/2} 2^{r\beta k}. 
\] 
Notice that $C_1(C,\beta,r)<\infty$ due to \eqref{eq:r_should_be_large}.
The lemma follows as we set $c=b2^{\gamma m}$ since there are $2^{(2+\alpha)m}$ points in $R_m$. \epf

The proof of the following lemma is similar to that of Lemma~\ref{lm:holder_from_moment1}, and we omit it.
\begin{lemma}\label{lm:holder_from_moment2}
 If $r/2-1>r\beta$ and $Z$ is a continuous process satisfying  
\begin{equation}
\label{eq:condition_on_increment_of_X}
 \E|Z(t')-Z(t)|^r\le C|t'-t|^{r/2},\quad t,t'\in[0,1],
\end{equation}
then
\[
 \Pp\{\Hc_\beta(Z)\ge b2^{\gamma m}\}\le C_1(C,\beta,r)\frac{2^{-r\gamma m}}{b^r}.
\]
\end{lemma}

\begin{lemma} Let 
\begin{equation}
\label{eq:condition_on_alpha}
2+\alpha<r\gamma.
\end{equation}
Suppose $\Gamma=(\Gamma,Z,U)$ is a monotone flow
such that processes $U$ satisfy the conditions of Lemma~\ref{lm:holder_from_moment1},
and process $Z$ satisfies the conditions of Lemma~\ref{lm:holder_from_moment2}.
Suppose there is a constant $\bar C$ such that
\begin{equation}
\label{eq:bound_on_supX}
\E \sup_{t\in[0,1]} Z(t)\le \bar C. 
\end{equation}
then there are $b=b(\alpha,\beta,\gamma,r,C,\bar C)$ and $m_0=m_0(\alpha,\beta,\gamma,r,C,\bar C)$ such that 
\[
\Pp\{\Gamma\in K\}>1-\eps.
\]
\end{lemma}
\bpf
\begin{multline}
 \Pp\{\Gamma\in K\}\le  \sum_{m=1}^\infty \Pp\left\{\sup_{(t_0,x)\in R_m} \Hc_\beta(U(x,t_0,\cdot)) > b2^{\gamma m}\right\}\\
+ \sum_{m=1}^{\infty} \Pp\left\{\Hc_\beta(Z)\ge b2^{\gamma m}\right\}\\+\sum_{m=m_0}^{\infty}\Pp\left\{ \sup_{t\in[0,1]} Z(t)>2^{\alpha m}\right\}.
\end{multline}
Lemma~\ref{lm:holder_from_moment1} and condition \eqref{eq:condition_on_alpha} imply that for sufficiently large $b$,
the first term in the r.h.s.\ is less than $\eps/3$. It also follows from Lemma~\ref{lm:holder_from_moment2} that for sufficiently large $b$, the second term in the r.h.s.\ is less than $\eps/3$. Finally,
we can use \eqref{eq:bound_on_supX} and Markov's inequality to choose $m_0$ so that the last term in the r.h.s.\ is less than $\eps/3$
completing the proof.
\epf

Since we want \eqref{eq:r_should_be_large} and \eqref{eq:condition_on_alpha} to be satisfied along with $\gamma<\beta$, we need to have
\[
2+\alpha<r\beta<\frac{r}{2}-1.
\]
This can be satisfied if we take $r=8$, $1/4<\gamma<\beta<3/8$, $0<\alpha<8\gamma-2$, and our goal is to show that all other conditions are satisfied with this choice of parameters for the sequence of random monotone flows $\Gamma_n=(\Gamma_n,Z_n,U_n)$ defined in~\eqref{eq:pre-limit-random-flows}. 

From now on we assume that $U_n$ is the special right-continuous trajectory representation
of $\Gamma_n$, see \eqref{eq:right_continuous_rep}.

\begin{lemma} For all $n$, the process $Z_n$ restricted to times $0,n^{-1},2n^{-1},3n^{-1},\ldots$ is a submartingale
with respect to its natural filtration. For any $r\in\N$, $\E Z_n^r(t)$ is uniformly bounded in  $t\in[0,1]$ and $n\in\N$.
\end{lemma}
\bpf
It is sufficient to consider the process $X_n$ counting vertices in the $n$-th generation of the discrete random tree and show that it
is a submartingale and that for every $r$ there is a constant $c_r$ such that for all $n\in\N$, 
\begin{equation}
\label{eq:moment_estimate_for_X^*}
\E X_n^r<c_rn^r.
\end{equation}
The latter follows from the following one: for every $r$ there is a polynomial $P_{r-1}$ of degree $r-1$
such that
\begin{equation}
\label{eq:iteration}
\E[X_{n+1}^r|X_n]\le  X_n^r+P_{r-1}(X_n). 
\end{equation}
Taking expectations of both sides of~\eqref{eq:iteration} and applying a straightforward induction procedure in $r$ and $n$,
we obtain the growth estimate~\eqref{eq:moment_estimate_for_X^*}.

To prove~\eqref{eq:iteration}, we define $P(r,k)$ to be the set of all partitions of $r$ into $k$ nonnegative integers:
\[
 P(r,k)=\{(j_1,\ldots,j_k):\ j_1\ge j_2\ge\ldots \ge j_k\ge 0,\ j_1+j_2+\ldots+j_k=r\},
\]
 and write
\begin{align}\notag
\E[X_{n+1}^r|X_n=k]
&=\frac{1}{k}\sum_{0\le i_1,\ldots,i_k\le D}(i_1+\ldots+i_k)^{r+1}p_{i_1}\ldots p_{i_k}
\\ \notag
&=\frac{1}{k}\sum_{0\le i_1,\ldots,i_k\le D}\sum_j\binom{r+1}{j_1,\ldots,j_k} i_1^{j_1}i_2^{j_2}\ldots i_k^{j_k}
p_{i_1}\ldots p_{i_k} 
\\ \notag
&=\frac{1}{k}\sum_{0\le i_1,\ldots,i_k\le D}\sum_{j\in P(r+1,k)}K_j(r)i_1^{j_1}i_2^{j_2}\ldots i_k^{j_k}
p_{i_1}\ldots p_{i_k}\\ 
&=\frac{1}{k}\sum_{j\in P(r+1,k)}K_j(r)B_{j_1}B_{j_{2}}\ldots B_{j_k}.
\label{eq:sum_in_terms_of_B}
\end{align}
Here $i_1,\ldots,i_k$ denote the number of children of each of $k$ vertices of generation $n$, and in the third line we used the symmetry
of the expression and grouped together monomials producing the same partition $j\in P(r+1,k)$ under the monotone rearrangement of their degrees (we agree that $0^0=1$). The numbers $B_m$ have been introduced in~\eqref{eq:definition_of_B_n}. The positive constants $K_{j}(r)$ satisfy 
\begin{equation}
\label{eq:sum_all_constants}
 \sum_{j\in P(r+1,k)}K_j(r)=k^{r+1}.
\end{equation}
Notice that if $k\ge r+1$ then 
\[
K_{1,1,\ldots,1,0,0\ldots 0}=k(k-1)\ldots (k-r)=k^{r+1}-\tilde P_r(k), 
\]
where $\tilde P_r(k)$ is a polynomial of degree at most $r$. In this case,
due to~\eqref{eq:sum_all_constants}, the sum of all other constants equals $\tilde P_r(k)$. Now \eqref{eq:iteration} follows
straight from~\eqref{eq:sum_in_terms_of_B}.

The fact that $(X_n)$ is a submartingale follows from
\[
 \E[X_{n+1}-X_{n}|\sigma(X_0,\ldots,X_{n})]=\mu,
\]
which was established in~\cite{Bakhtin-RSA}. This identity can also be derived from the computation in the proof of 
Lemma~\ref{lm:conditions_are_satisfied}.
\epf

\begin{lemma}\label{lm:conditions_are_satisfied} The conditions of Lemma~\ref{lm:holder_from_moment1} are satisfied for
$U_n$ uniformly in~$n$, if we choose $r=8$, $1/4<\gamma<\beta<3/8$, $0<\alpha<8\gamma-2$. 
\end{lemma}
\bpf We have to prove that~\eqref{eq:moment_estimate} is satisfied for $U_n(x,t_0,\cdot)$ uniformly in $n,x,t_0,t,t'$.
Suppose first that $(t_0,x)$ is a grid point, i.e., $t_0 n$ and $x\mu n$ are integers. Denote by $V_j$ the size of progeny in generation~$j$ generated by first $x\mu n+1$ vertices in $t_0n$-th generation. Then $U_n(x,t_0,j/n)=(V_j-1)/(\mu n)$. 
Let us compute
\begin{align*}
&\E\Biggl[\frac{V_{m+1}}{\mu n} \Biggr| \frac{V_{m}}{\mu n}=\frac{l}{\mu n},\ \frac{X_m}{\mu n}=\frac{k}{\mu n}\Biggr]\\
&= \frac{1}{\mu nk}\sum_{0\le i_1,\ldots,i_k\le D}(i_1+\ldots+i_l)(i_1+\ldots+i_k)p_{i_1}\ldots p_{i_k}\\
&=\frac{1}{\mu nk}\left[l\sum_{i_1}i_1^2 p_{i_1}\ldots p_{i_k}
+l(k-1)\sum_{i_1\ne i_2}i_1i_2p_{i_1}\ldots p_{i_k} \right]\\
&=\frac{1}{\mu nk}\left[l B_2+l(k-1)\right]\\
&=\frac{l}{\mu n}\left(1+\frac{\mu}{k}\right).
\end{align*}
Therefore,
\begin{equation}
\E\Biggl[\frac{V_{m+1}}{\mu n} - \frac{V_{m}}{\mu n} \Biggr| \frac{V_{m}}{\mu n}=\frac{l}{\mu n},\ \frac{X_m}{\mu n}=\frac{k}{\mu n}\Biggr]=\frac{1}{n}\frac{l}{k}.
\label{eq:increment_expect}
\end{equation}
Since $0\le l/k\le 1$, the process $\bar A_{n}(t)$ interpolating between the values of
\[
 A_{n}(t)=\frac{1}{\mu n}\sum_{j=0}^{[nt]-1}\frac{V_{j}}{X_{j}}
\]
is a 1-Lipschitz process, and
\[
 B_n(t)=U_n(x,t_0,t)-\bar A_n(t)
\]
restricted to $t=mn^{-1}$, $m=0,1,\ldots,n$, is a martingale. Since $\bar A$ is 1-Lipschitz, we now need to prove that $B_n$ satisfies a moment estimate analogous to~\eqref{eq:moment_estimate}
with constant in the r.h.s.\ independent of $n$. If $t=mn^{-1}$ and $t'=m'n^{-1}$, then due to Burkholder's inequality 
(see \cite[Section VII.3]{Shiryaev:MR1368405}) and Minkowski's inequality (see \cite[II.6]{Shiryaev:MR1368405}), there is an absolute constant $C_2$ such that
\begin{align*}
 \E |B_n(t')-B_n(t)|^8&\le C_2 \E \left(\sum_{j=m}^{m'-1}\left(B_n((j+1)n^{-1})-B_n(jn^{-1})\right)^2\right)^4
\\&\le  C_2 \left\|\sum_{j=m}^{m'-1}\left(B_n((j+1)n^{-1})-B_n(jn^{-1})\right)^2\right\|_4^4\\
&\le C_2 \left(\sum_{j=m}^{m'-1}\left\|\left(B_n((j+1)n^{-1})-B_n(jn^{-1})\right)^2\right\|_4\right)^4\\
&\le C_2(m'-m)^4\max_{j<m'} \left\|\left(B_n((j+1)n^{-1})-B_n(jn^{-1})\right)^2\right\|_4^4\\
&\le C_2 n^4 (t'-t)^4\max_{j<n} \E\left|B_n((j+1)n^{-1})-B_n(jn^{-1})\right|^8.
\end{align*}
To estimate the r.h.s., we use the following statement: 
\begin{lemma}\label{lm:one-step_moment_estimate} There is an absolute constant $C_3$ such that
 \[\E\left|B_n((j+1)n^{-1})-B_n(jn^{-1}))\right|^8\le \frac{C_3}{n^4}.\]
\end{lemma}
The desired moment estimate follows directly:
\[
 \E |B_n(t')-B_n(t)|^8\le C_2C_3(t'-t)^4. 
\]

\bpf[Proof of Lemma~\ref{lm:one-step_moment_estimate}]
\[
\left|B_n\left(\frac{j+1}{n}\right)-B_n\left(\frac{j}{n}\right)\right|\le \frac{1}{\mu n}+
\left|\frac{V_{j+1}}{\mu n}-\frac{V_{j}}{\mu n}\right|,
\]
so that, by convexity,
\begin{align*}
\E\left|B_n\left(\frac{j+1}{n}\right)-B_n\left(\frac{j}{n}\right)\right|^8&\le 2^7\frac{1}{(\mu n)^8}+
2^7\E\left|\frac{V_{j+1}}{\mu n}-\frac{V_{j}}{\mu n}\right|^8,
\end{align*}
and the lemma follows directly from the next result:
\begin{lemma}\label{lm:moment_estimate_increment}
For any even $r\ge 2$, there is a constant $c'_r$ such that for all $n$ and all $j$,
\begin{equation}
 \E \left(\frac{V_{j+1}}{n}-\frac{V_j}{n}\right)^r\le  c'_r \frac{j^{r/2}}{n^{r}}.
\label{eq:moment_estimate2}
\end{equation}
\end{lemma}
\bpf Let us estimate $\E(V_{j+1}-V_j)^r$.
\begin{align*}
 \E [(V_{j+1}&-V_j)^r|\ V_j=l,X_j=k]\\
=&\frac{1}{k}\sum_{i_1,\ldots,i_k}(i_1+\ldots+i_l-l)^r(i_1+\ldots+i_k-k)p_{i_1}\ldots p_{i_k}\\
&+\sum_{i_1,\ldots,i_k}(i_1+\ldots+i_l-l)^rp_{i_1}\ldots p_{i_k}\\
=&\frac{1}{k}\sum_{i_1,\ldots,i_k}((i_1-1)+\ldots+(i_l-1))^r((i_1-1)+\ldots+(i_k-1))p_{i_1}\ldots p_{i_k}\\
&+\sum_{i_1,\ldots,i_k}((i_1-1)+\ldots+(i_l-1))^rp_{i_1}\ldots p_{i_k}
\\=&\frac{1}{k}\sum_{i_1,\ldots,i_k}((i_1-1)+\ldots+(i_l-1))^{r+1}p_{i_1}\ldots p_{i_k}\\
&+\frac{1}{k}\sum_{i_1,\ldots,i_k}((i_1-1)+\ldots+(i_l-1))^r((i_{l+1}-1)+\ldots+(i_k-1))p_{i_1}\ldots p_{i_k}\\
&+\sum_{i_1,\ldots,i_k}((i_1-1)+\ldots+(i_l-1))^rp_{i_1}\ldots p_{i_k}
\\ =&S_1+S_2+S_3
\end{align*}
We begin with $S_3$. Since
\[
 \sum_{i_{l+1},\ldots,i_k}p_{i_{l+1}}\ldots p_{i_k}=1,
\]

\begin{align*}
 S_3&=\sum_{j\in P(r,l)}\sum_{i_1,\ldots,i_l}
(i_1-1)^{j_1}\ldots(i_l-1)^{j_l}p_1\ldots p_l\\
&=\sum_{j\in P(r,l)}K_{j}\bar B_{j_1}\ldots \bar B_{j_l},
\end{align*}
where all monomials with similar degree $j$ are grouped together, and
\[
 \bar B_{j}=\sum_{i}(i-1)^jp_i.
\]
Notice that $\bar B_{1}=0$, so that if one of the terms in $j_1+\ldots+j_l=r$ equals 1, the contribution from
that term to $S_3$ is 0. In other words, only partitions of $r$ containing no $1$'s contribute to $S_3$. Each of these
partitions contains at most $r/2$ nonzero components, and therefore the associated coefficient $K_{j}$ 
(the number of monomials associated to this partition)
is
bounded by a degree $r/2$ polynomial in $l$.
The same analysis shows that $S_2=0$, and $S_1$ is bounded by a polynomial of degree $r/2$. We conclude
that there is a constant $C'_r$ (independent of $k$) such that
\begin{equation}
\label{eq:one_step_moment_estimate_cond_on_l}
 \E [(V_{j+1}-V_j)^r|\ V_j=l,X_j=k]\le C'_rl^{r/2}.
\end{equation}
The estimate~\eqref{eq:moment_estimate_for_X^*} implies now that
\[
 \E (V_{j+1}-V_j)^r\le  C'_r\E X_j^{r/2}\le c'_r j^{r/2},
\]
so that~\eqref{eq:moment_estimate2} holds true.
\epf

We assume now that the desired moment estimate for $(t_0,x)$ on the grid holds true with a constant $C$.
We need is to estimate 
$\E|U_n(x,t_0,t')-U_n(x,t_0,t)|^8$ for arbitrary $n,x,t_0,t,t'$.

Let $j=[tn]$, $j'=[t'n]$ and $N=j'-j$.  We shall assume that
$t_0n < j$; the proof can be easily modified to treat the opposite situation.

We are going to consider three cases:  (i) $N=0$, (ii) $N=1$, (iii) $N\ge 2$.

In case (i), the evolution of $U_n(x,t_0,t)$ is linear in $t\in[jn^{-1},(j+1)n^{-1}]$, so that, comparing
it to the evolution along the regular edges of the tree imbedding, we see that 
\begin{align*}
 &\E \left(\frac{U_n(x,t_0,t')-U_n(x,t_0,t)}{t'-t}\right)^8  =\E\left(\frac{ U_n(x,t_0,\frac{j+1}{n})-U_n(x,t_0,\frac{j}{n})}{1/n}\right)^8\\
\le& n^8 \sum_{l}\E\left[ (U_n(x,t_0,\tst\frac{j+1}{n})-U_n(x,t_0,\frac{j}{n}))^8|\ U_n(x,t_0,\frac{j}{n})\in(\frac{l}{\mu n},
\textstyle \frac{l+1}{\mu n}]\right],  \\
\le& \sum_{l}\left(\E [(V_{j+1}-V_j)^8| V_{j}=l+1]\Pp\left\{\tst U_n(x,t_0,\frac{j}{n})\in(\frac{l}{\mu n},\frac{l+1}{\mu n}]\right\}\right.\\ &+
\left.\E [(V_{j+1}-V_j)^8| V_{j}=l]\Pp\left\{\tst U_n(x,t_0,\frac{j}{n})\in(\frac{l}{\mu n},\frac{l+1}{\mu n}]\right\}\right). 
\end{align*}
Inequality \eqref{eq:one_step_moment_estimate_cond_on_l} implies now that we can continue the estimate as
\begin{align}\notag
\E(U_n(x,t_0,t')&-U_n(x,t_0,t))^8
\\&\le 2(t'-t)^8 C'_8\sum_{l}(l+1)^4 \Pp\left\{\tst U_n(x,t_0,\frac{j}{n})\in(\frac{l}{\mu n},\frac{l+1}{\mu n}]\right\}
\notag
\\&\le 2(t'-t)^8\mu^4 n^4 C'_8\, \E\, \tst (U_n(x,t_0,\frac{j}{n})+\frac{1}{\mu n})^4\notag
\\&\le 2(t'-t)^4 \mu^4 C'_8 C_4^*,\label{eq:case_i}
\end{align} 
where 
\[
C_4^*=\sup_{n\in\N}\sup_{t\in[0,1]}\E \left(Z_n(t)+\frac{1}{\mu n}\right)^4, 
\]
and in the last line we used the fact that $t'-t\le n^{-1}$.
In case (ii),
\begin{multline*}
 \E\left(\frac{U_n(x,t_0,t')-U_n(x,t_0,t)}{t'-t}\right)^8\\ \le \E\left(\frac{U_n(x,t_0,\frac{j+1}{n})-U_n(x,t_0,t)}{\frac{j+1}{n}-t}\right)^8
+\E\left(\frac{U_n(x,t_0,t')-U_n(x,t_0,\frac{j+1}{n})}{t'-\frac{j+1}{n}}\right)^8
\\ \le \E\left(\frac{U_n(x,t_0,\frac{j+1}{n})-U_n(x,t_0,\frac{j}{n})}{1/n}\right)^8
+\E\left(\frac{U_n(x,t_0,\frac{j+2}{n})-U_n(x,t_0,\frac{j+1}{n})}{1/n}\right)^8, 
\end{multline*}
and for each of the terms in the r.h.s.\ we can proceed as in case (i) and obtain an estimate analogous to~\eqref{eq:case_i}. 

Case (iii). Let us estimate $\E(U_n(x,t_0,t')-U_n(x,t_0,t))_+^8$ first.
Let $\xi=[U(x,t_0,t)\mu n]$. Suppose $U(x,t_0,t')\ge U(x,t_0,t)$. Then the monotonicity of the flow implies that at least one of the following two
conditions holds true:
\begin{enumerate}
 \item 
$
U( \frac{\xi}{\mu n},\frac{j}{n},\frac{j+1}{n})- U(x,t_0,t)\ge  \frac{1}{2}(U(x,t_0,t')-U(x,t_0,t)). 
$
\item For $\eta=[U(x,t_0,\frac{j+1}{n})\mu n]+1$,
\[\tst
  U(\frac{\eta}{\mu n},\frac{j+1}{n},\frac{j'-1}{n})\wedge U(\frac{\eta}{\mu n},\frac{j+1}{n},\frac{j'}{n}) - \frac{\eta}{\mu n}\ge \frac{1}{2}(U(x,t_0,t')-U(x,t_0,t))-\frac{1}{\mu n}.
\]
\end{enumerate} 
Consequently,
\[
 (U(x,t_0,t')-U(x,t_0,t))_+\le 2(R_1\wedge R_2 \wedge R_3), 
\]
where
\begin{align*}
R_1&=U\left(\frac{\xi}{\mu n},\frac{j}{n},\frac{j+1}{n}\right)-U(x,t_0,t),\\
R_2&=U\left(\frac{\eta}{\mu n},\frac{j+1}{n},\frac{j'-1}{n}\right)- \frac{\eta}{\mu n}+\frac{1}{n}\\
R_3&=U\left(\frac{\eta}{\mu n},\frac{j+1}{n},\frac{j'}{n}\right) - \frac{\eta}{\mu n}+\frac{1}{n}.
\end{align*}
Therefore,
\[
 \E(U(x,t_0,t')-x)_+^8\le 2^8(\E R_1^8+ \E R_2^8 +\E R_3^8). 
\]
Notice that the random variables $R_1,R_2,R_3$ are defined in terms of the flow trajectories emitted from lattice points.
Although these lattice points defined through $\xi$ and $\eta$ are random, they are measurable with respect to the filtration
of the flow, and therefore the above estimates imply that there is a constant $C$ such that
\[
 \E(U(x,t_0,t')-U(x,t_0,t))_+^8\le C\left(\frac{j'-j}{n}\right)^4\le 2^4 C(t'-t)^4.
\]
We can estimate $\E(U(x,t_0,t')-U(x,t_0,t))_-^8$ similarly,  which completes the proof.

\subsection{Characterization of limit points}\label{sbsc:characterization}

We begin with two auxiliary lemmas.
\begin{lemma}\label{lm:small_v_stays_small}  Let $m_0\in\N$.
For each $n\in\N$, take a random number $V(n)$ of vertices in generation $m_0n$ measurable
with respect to the history of the tree up to generation $m_0n$. 
For $m\ge m_0n$ we denote by
$V_m^{(n)}$ the total size of the progeny of these vertices in generation $m$. (In particular,
$V_{m_0n}^{(n)}=V(n)$.)
If
\[
\Pp\{V{(n)}/n\ge n^{-\gamma}\}\to 0,\quad n\to\infty, 
\]
for some constant $\gamma\in(0,1)$, independent of $n$, then, for any $m_1>m_0$,
\[
 \sup_{m_0n\le m\le m_1n} \frac{V_m^{(n)}}{n}\stackrel{\Pp}{\to} 0,\quad n\to\infty.
\]
\end{lemma}
\bpf Let
$X_m, m\ge 0$ is the process of total population sizes.
If we show that 
\begin{equation}
\label{eq:ratio_V_X_conv_to0}
\sup_{m_0n\le m\le m_1n} \frac{V^{(n)}_m}{X_m}\stackrel{\Pp}{\to} 0,\quad n\to\infty, 
\end{equation}
then the lemma will follow since we can write
\[
\sup_{m_0n\le m\le m_1n} \frac{V^{(n)}_m}{n}= \sup_{m_0n\le m\le m_1n} \frac{V^{(n)}_m}{X_m} \cdot \sup_{m_0n\le m\le m_1n} \frac{X_m}{n}, 
\]
and use the convergence in distribution of $\sup_{m_0n\le m\le m_1n} \frac{X_{m}}{n}$ to the maximum of a diffusion process. 

To prove~\eqref{eq:ratio_V_X_conv_to0}, we define a sequence of events \[A_n= \{X_{m_0n}\ge n^{1-\gamma/2};\ V(n)\le n^{1-\gamma}\},\]
and for any $\eps>0$ write
\[
\Pp\left\{\sup_{m_0n\le m\le m_1n} \frac{V_{m}}{X_{m}}>\eps \right\}\le \Pp(A_n^c)+ \Pp\left(\left\{\sup_{m_0n\le m\le m_1n} \frac{V_{m}}{X_{m}}>\eps \right\}\cap A_n\right).
\]

Clearly, $\Pp(A_n^c)\to 0$ as $n\to\infty$ due to the assumptions of the lemma and the fact that the distribution of $X_n/n$ converges weakly
to a distribution with no atom at $0$. 

For the second term, we notice that $\left(\frac{V^{(n)}_m}{X_m}\right)$ is a nonnegative bounded supermartingale w.r.t.\ filtration~$(\Fc_m^{(n)})$, where for each $m=m_0n,\ldots,m_1n$,  $\Fc_m^{(n)}$ is the sigma-algebra generated by $X_l,V_l,l=m_0n,\ldots,m$. 
In fact,

\begin{align*}
 \E &\left[\frac{V^{(n)}_{m+1}}{X_{m+1}}|\ V^{(n)}_m=l, X_m=k\right]
\\&=\frac{1}{k}\sum_{\substack{0\le i_1,\ldots,i_k\le D\\ i_1+\ldots+ i_k\ne0}}\frac{i_1+\ldots+ i_l}{i_1+\ldots+ i_k}
(i_1+\ldots+ i_k)p_{i_1}\ldots p_{i_k}
\\&=\frac{1}{k}\sum_{\substack{0\le i_1,\ldots,i_k\le D}}(i_1+\ldots+ i_l)p_{i_1}\ldots p_{i_k}
-\frac{1}{k}{p_0}^k
\\&=\frac{l}{k}-\frac{1}{k}{p_0}^k
\\&\le \frac{l}{k}.
\end{align*}
Doob's maximal inequality for nonnegative supermartingales implies 
\[
 \Pp\left(\left\{\sup_{m_0n\le m\le m_1n} \frac{V_{m}}{X_{m}}>\eps \right\}\cap A_n\right)\le \frac{n^{1-\gamma}/ n^{1-\gamma/2}}{\eps}\to0,\quad n\to\infty,
\]
and~\eqref{eq:ratio_V_X_conv_to0} is proven. 
\epf

\begin{lemma} \label{lm:small_perturb_if_start_not_on_grid}Let $\gamma\in(0,1/2)$.
 For any $x>0,t>0$, there is a constant $C$ such that  for all $n\in\N$,
\[
 \Pp\left\{\tst |U_n(x,t,\frac{[nt]+1}{n})- \frac{[x\mu n]}{\mu n}|>n^{-\gamma}\right\}< Cn^{-(1-2\gamma)}
\]
\end{lemma}
\bpf It is easily seen that for large $n$, due to the monotonicity of the flow, if 
\[
 \left|U_n\left(x,t,\frac{[nt]+1}{n}\right)- \frac{[x\mu n]}{\mu n}\right|>n^{-\gamma},
\]
then 
\[
  U_n\left(\frac{[x\mu n]}{\mu n},t,\frac {[nt]+1}{n}\right)- \frac{[x\mu n]}{\mu n}>0,
\]
and, moreover,
\[
\left| U_n\left(\frac{[x\mu n]}{\mu n},t,\frac{[nt]+1}{n}\right)- \frac{[x\mu n]}{\mu n}\right| > n^{-\gamma}-\frac{2}{n}.
\]
Also, for large values of $n$, if
\[
U_n\left(x,t,\frac{[nt]+1}{n}\right)- \frac{[x\mu n]}{\mu n} < - n^{-\gamma}. 
\]
then 
\[
U_n\left(\frac{[x\mu n]}{\mu n},t,\frac{[nt]+1}{n}\right)- \frac{[x\mu n]}{\mu n}<0 
\]
and, moreover,
\[
\left|U_n\left(\frac{[x\mu n]}{\mu n},t,\frac{[nt]+1}{n}\right)- \frac{[x\mu n]}{\mu n} \right| > n^{-\gamma}-\frac{2}{n}.
\]
Therefore,
\begin{multline*}
 \Pp\left\{\tst |U_n\left(x,t,\frac{[nt]+1}{n}\right)- \frac{[x\mu n]}{\mu n}| > n^{-\gamma}\right\}
\\ \le \Pp\left\{\tst \left|U_n\left(\frac{[x\mu n]}{\mu n},t,\frac{[nt]+1}{n}\right)-\frac {[x\mu n]}{\mu n} \right| > n^{-\gamma}-\frac{2}{n}\right\},
\end{multline*}
The moment estimate~\eqref{eq:one_step_moment_estimate_cond_on_l} from Lemma~\ref{lm:moment_estimate_increment} implies that 
\begin{align*}
 \E \left|U_n\left(\frac{[x\mu n]}{\mu n},t,\frac{[nt]+1}{n}\right)-\frac {[x\mu n]}{\mu n} \right|^2&
\le C'_2\frac{[x\mu n]}{\mu^2 n^2}
\le C'_2\frac{x}{n},	 
\end{align*}
and the desired estimate follows from Markov's inequality.
\epf

\bigskip

Now we proceed to prove that any limiting point for the sequence of distributions of monotone flows $\Gamma_n$
has to coincide with the monotone flow solving the SPDE~\eqref{eq:spde-1-integral} as discussed in Section~\ref{sec:spde_solution}.

Consider $(x_{0,1},t_0)=(0,0)$, a sequence of times $0<t_1<t_2<\ldots<t_k$, and for each $i=1,\ldots,k$, a sequence of
nonnegative numbers \[x_{i1}<\ldots<x_{il(i)}.\] 
For each $i=0,\ldots,k$ and every $j=1,\ldots, l(i)$ we define a family of processes
\[
 U^{i,j}_n(t)=\begin{cases}
         x_{ij},& t<\frac{[n t_i ]}{n},\\
         U_n\left(\frac{[x_{ij}\mu n]}{\mu n},\frac{[nt_i]+1}{n},t\right), & t\ge\frac{[nt_i]}{n}.
        \end{cases}
\]

Due to Lemmas \ref{lm:small_v_stays_small} and \ref{lm:small_perturb_if_start_not_on_grid}, it is sufficient to show that as $n\to\infty$, these processes jointly converge in distribution in sup-norm to a nonnegative diffusion process 
$(U^{ij}_\infty)$
with drift
\[
b^{ij}(y,t)=  \frac{y_{ij}}{y_{01}}\ONE_{\{t>t_i;\ 0\le y_{ij}\le y_{01}\}},\quad i=0,\ldots,k,\quad j=1,\ldots, l(i), 
\]
and diffusion matrix
\begin{multline*}
 a^{i_1j_1, i_2j_2}(y,t)= (y_{i_1j_1}\wedge y_{i_2j_2})\ONE_{\{t>t_{i_1}\vee t_{i_2};\ 0\le y_{i_1j_1},y_{i_2j_2}\le y_{01}\}},\\
             i_1,i_2=0,\ldots,k,\quad j_1=1,\ldots, l(i_1),\quad j_2=1,\ldots, l(i_2).
\end{multline*}
For any $i$ and any $y$, these coefficients are constant on $t\in[t_i,t_{i+1})$, they are continuous and bounded in $y$,  and define a well-posed martingale problem.

Let 
\[
\bar U^{ij}_n(t)= U^{ij}_n([nt]/n) 
\]
If $t>t_i$ and $U^{ij}_n(\frac{[nt]}{n})\le U^{01}_n(\frac{[nt]}{n})$,
\[
 \E \left[U^{ij}_n\left(\frac{[nt]}{n}+\frac{1}{n}\right)|\ \Fc_{[nt]/n}\right]= U^{ij}_n\left(\frac{[nt]}{n}\right)+\frac{1}{n}\cdot\frac{U^{ij}_n(\frac{[nt]}{n})}{U^{01}_n(\frac{m}{n})}.
\]
So, we define
\begin{equation}
 B^{ij}_n(t)=\frac{1}{n}\sum_{m=0}^{[nt]}\ONE_{\left\{m\ge nt_i; U^{ij}_n(\frac{[nt]}{n})\le U^{01}_n(\frac{[nt]}{n})\right\}}\frac{U^{ij}_n(\frac{m}{n})}{U^{01}_n(\frac{m}{n})}
\label{eq:compensator}
\end{equation}
and
\[
 M^{ij}_n(t)=\bar U^{ij}_n(t)-B^{ij}_n(t-\frac{1}{n}).
\]
Next, a simple calculation based on the martingale property of $M^{ij}_n$ and~\eqref{eq:compensator} shows that
if $t>t_i$, $U^{i_1j_1}_n(\frac{[nt]}{n})\le U^{01}_n(\frac{[nt]}{n})$, and  $U^{i_2j_2}_n(\frac{[nt]}{n})\le U^{01}_n(\frac{[nt]}{n})$
then
\begin{align}\notag
\E&\left[M^{i_1j_1}_n\left(\frac{[nt]}{n}+\frac{1}{n}\right)M^{i_2j_2}_n\left(\frac{[nt]}{n}+\frac{1}{n}\right)
-M^{i_1j_1}_n\left(\frac{[nt]}{n}\right)M^{i_2j_2}_n\left(\frac{[nt]}{n}\right)|\ \Fc_{[nt]/n}\right]\\
=&\E\left[\bar U^{i_1j_1}_n\left(\frac{[nt]}{n}+\frac{1}{n}\right)\bar U^{i_2j_2}_n\left(\frac{[nt]}{n}+\frac{1}{n}\right) |\ \Fc_{[nt]/n}\right]\label{eq:compute_predictable_quadratic_covariation}
\\ \notag &-
\bar U^{i_1j_1}_n\left(\frac{[nt]}{n}\right)\left(1+\frac{1}{n}\cdot\frac{B_2-1}{U^{01}_n\left(\frac{[nt]}{n}\right)}\right)
\bar U^{i_2j_2}_n\left(\frac{[nt]}{n}\right)\left(1+\frac{1}{n}\cdot\frac{B_2-1}{U^{01}_n\left(\frac{[nt]}{n}\right)}\right).
\end{align}

Processes $U^{i_1j_1}_n$ and $U^{i_2j_2}_n$ describe the evolution of sizes of two subpopulations, i.e., they are rescaled
versions of vertex-counting discrete processes $V^{(1)}$ and $V^{(2)}$.  
To compute the first term in the r.h.s.\ of \eqref{eq:compute_predictable_quadratic_covariation}, we take $0\le l_1\le l_2\le k$ and write
\begin{align*}
\E&\left[V^{(1)}_{m+1}V^{(2)}_{m+1}|\ V^{(1)}_{m}=l_1, V^{(2)}_{m}=l_2, X_m=k\right]
\\&=\frac{1}{k}\sum_{i_1,\ldots,i_k}
(i_1+\ldots+i_{l_1})(i_1+\ldots+i_{l_2})(i_1+\ldots+i_k)p_{i_1}\ldots p_{i_k}
\\&=\frac{1}{k}(l_1B_3+(l_1(k-1)+l_1(l_2-1)+l_1(l_2-1))B_2+l_1(l_2-1)(k-2)).
\end{align*}
Denoting 
\begin{align*}
l_1&=\mu n\left(\bar U^{i_1j_1}_n\left(\frac{[nt]}{n}\right)\wedge \bar U^{i_2j_2}_n\left(\frac{[nt]}{n}\right)\right),\\
l_2&=\mu n\left(\bar U^{i_1j_1}_n\left(\frac{[nt]}{n}\right)\vee \bar U^{i_2j_2}_n\left(\frac{[nt]}{n}\right)\right),\\
k&=\mu n\bar U^{01}_n\left(\frac{[nt]}{n}\right),
\end{align*}
and making a simple calculation,
we can rewrite the r.h.s.\ of~\eqref{eq:compute_predictable_quadratic_covariation} as
\begin{align*}
\frac{l_1}{n}\cdot\left(\frac{1}{n}+\frac{R(k,l_2)}{kn}\right),
\end{align*}
where $R(\cdot,\cdot)$ is a uniformly bounded function. We denote
\begin{align}\label{eq:quadratic_predictable_def}
 A_n^{i_1j_1i_2j_2}(t)=\frac{1}{n}\sum_{m=0}^{[nt]}&\ONE_{\left\{m\ge n(t_{i_1}\vee t_{i_2});
U^{i_1j_1}_n(\frac{[nt]}{n}), U^{i_2j_2}_n(\frac{[nt]}{n})\le U^{01}_n(\frac{[nt]}{n})
\right\}}
\times \\ &\times  \left(\bar U^{i_1j_1}_n\left(\frac{m}{n}\right)\wedge \bar U^{i_2j_2}_n\left(\frac{m}{n}\right)\right)\times
\\ &\times \left(1+
\frac{
R(
\bar U_n^{01}(\frac{m}{n}), n(\bar U^{i_1j_1}_n(\frac{m}{n})\vee \bar U^{i_2j_2}_n(\frac{m}{n}))
)
}{n\bar U_n^{01}(\frac{m}{n})}\right).\notag
\end{align}
The above calculation shows that
\[
  M^{i_1j_1}_n(t)M^{i_1j_1}_n(t)-B_n^{i_1j_1i_2j_2}(t)
\]
is a martingale. 

We are going to use Theorem 4.1 from \cite[Chapter 7]{Ethier-Kurtz:MR838085}.
Identities~\eqref{eq:compensator} and
~\eqref{eq:quadratic_predictable_def} imply two main conditions of that theorem:
\[
 \sup_{0\le t\le 1}\left|B_n^{ij}(t)-\int_{0}^tb_n^{ij}(\bar U(s),s)ds\right|\stackrel{\Pp}{\to} 0,
\]
and
\[
\sup_{0\le t\le 1}\left|A_n^{i_1j_1i_2j_2}(t)-\int_{0}^t a_n^{i_1j_1i_2j_2}(\bar U(s),s)ds\right|\stackrel{\Pp}{\to} 0.
\]
The other conditions of Theorem 4.1\cite{Ethier-Kurtz:MR838085} are concerned with the jumps of processes $A$ and $B$, they can be easily
proven using the moment estimates above. The theorem allows us 
 to conclude that
the process $\bar U_n$ converges weakly in the Skorokhod topology on $D[0,1]$ to
the continuous diffusion process with drift and diffusion given by $b$ and $a$ respectively. In fact, Theorem~4.1 is given in
\cite{Ethier-Kurtz:MR838085} for time-independent coefficients, and here we invoke its straightforward time-dependent generalization.

Convergence of $\bar U_n$ in distribution in Skorokhod topology to a continuous process  implies convergence in distribution in sup norm to the same process. Moreover we can conclude that the $U_n$, the linear interpolated version of $\bar U_n$ also converges in distribution in sup norm,
which concludes the demonstration of the theorem. 

\section{Connection to superprocesses}\label{sec:superprocesses}

There is an important connection of our results to the theory of superprocesses. 
Superprocesses are measure-valued stochastic processes
describing the evolution of populations of branching and migrating particles, see e.g.~\cite{Dawson:MR1242575}. The limiting SPDE that we have constructed is similar to
the genealogy in the Dawson--Watanabe superprocess with no motion conditioned on nonextinction, see~\cite{Evans:MR1249698},\cite{Evans-Perkins:MR1088825},
\cite{Donnelly-Kurtz:MR1728556}. 

Our approach is more
geometric then the superprocess point of view. For the superprocess corresponding to our situation, the continual mass momentarily organizes
itself into a finite random number of atoms of positive mass (corresponding to discontinuities of the monotone maps in our approach). The mass of these atoms evolves in time analogously to 
equations~(\ref{eq:U-time-noise}), but our approach helps to understand what happens inside the atoms by unfolding the details of the genealogy giving
rise to a contentful
cocycle property whereas in superprocesses the cocycle action essentially reduces to evolving the masses of atoms. The monotone flow we construct
can be initiated with zero mass and just one particle whereas in superprocesses one has to start with positive mass.
The related historical superprocesses do describe the genealogy, but they require nontrivial particle motion and ignore the 
ordering and monotonicity issues that are central to our approach.

Notice that the continuous monotone flow we construct, describes a foliation of the random set
$\{ (t,x):\ 0\le x\le Z(t)\}$ into diffusion trajectories. The geometry of this stochastic foliation is suprisingly complicated.
Due to the presence of shocks in monotone maps (corresponding to the atoms of superprocesses), the stochastic 
foliation cannot be obtained, say, as a continuous image of the foliation of a rectangle in ``horizontal'' segments, 
and the flow is very far from a flow of diffeomorphisms.

We hope that our results are interesting from the point of view of graph theory since they describe what a typical large tree looks like.
In fact, combining the results of this paper with those of Section~\ref {sec:Gibbs}, we conclude that a typical embedding of a large ordered rooted tree in
the plane if rescaled appropriately looks like a stochastic foliation described by SPDE~(\ref{eq:spde-with-X*}). It would be interesting to obtain
rigorously a direct convergence result that would not involve the intermediate infinite discrete tree. However, currently this kind of result
is not available.

\begin{figure}
\centering
\includegraphics[height=4.4cm]{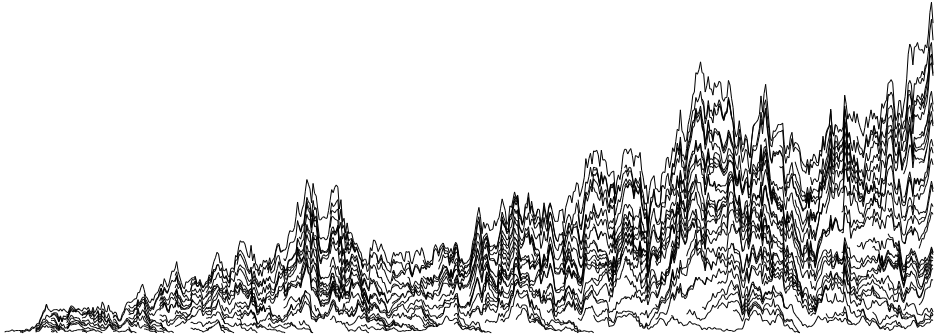}
\caption{Stochastic foliation constructed for 600 generations of a tree.}
\label{fig:stochastic_foliation-2}
\end{figure}

Figure~\ref{fig:stochastic_foliation-2} shows a realization of a pre-limit monotone flow for a large random tree. 
Every tenth generation is split into about ten
subpopulations, their progenies are tracked and shown on the figure.

\bibliographystyle{alpha}
\bibliography{treelimit}

\begin{thebibliography}{Daw93}

\bibitem[AP98]{Aldous-Pitman:MR1641670}
David Aldous and Jim Pitman.
\newblock Tree-valued {M}arkov chains derived from {G}alton-{W}atson processes.
\newblock {\em Ann. Inst. H. Poincar\'e Probab. Statist.}, 34(5):637--686,
  1998.

\bibitem[Bak]{Bakhtin-RSA}
Yuri Bakhtin.
\newblock Thermodynamic limit for large random trees.
\newblock {\em Accepted at Random Structures Algorithms}.
\newblock Also available at http://arxiv.org/abs/0809.2974.

\bibitem[BH08]{Bakhtin-Heitsch-1:MR2415118}
Yuri Bakhtin and Christine Heitsch.
\newblock Large deviations for random trees.
\newblock {\em J. Stat. Phys.}, 132(3):551--560, 2008.

\bibitem[BH09]{Bakhtin-Heitsch-2}
Yuri Bakhtin and Christine Heitsch.
\newblock Large deviations for random trees and the branching of {RNA}
  secondary structures.
\newblock {\em Bull Math Biol}, 71(1):84--106, 2009.

\bibitem[Daw93]{Dawson:MR1242575}
Donald~A. Dawson.
\newblock Measure-valued {M}arkov processes.
\newblock In {\em \'{E}cole d'\'{E}t\'e de {P}robabilit\'es de {S}aint-{F}lour
  {XXI}---1991}, volume 1541 of {\em Lecture Notes in Math.}, pages 1--260.
  Springer, Berlin, 1993.

\bibitem[DK99]{Donnelly-Kurtz:MR1728556}
Peter Donnelly and Thomas~G. Kurtz.
\newblock Genealogical processes for {F}leming-{V}iot models with selection and
  recombination.
\newblock {\em Ann. Appl. Probab.}, 9(4):1091--1148, 1999.

\bibitem[EK86]{Ethier-Kurtz:MR838085}
Stewart~N. Ethier and Thomas~G. Kurtz.
\newblock {\em Markov processes. Characterization and convergence}.
\newblock Wiley Series in Probability and Mathematical Statistics: Probability
  and Mathematical Statistics. John Wiley \& Sons Inc., New York, 1986.

\bibitem[EP90]{Evans-Perkins:MR1088825}
Steven~N. Evans and Edwin Perkins.
\newblock Measure-valued {M}arkov branching processes conditioned on
  nonextinction.
\newblock {\em Israel J. Math.}, 71(3):329--337, 1990.

\bibitem[Eva93]{Evans:MR1249698}
Steven~N. Evans.
\newblock Two representations of a conditioned superprocess.
\newblock {\em Proc. Roy. Soc. Edinburgh Sect. A}, 123(5):959--971, 1993.

\bibitem[Ken75]{Kennedy:MR0386042}
Douglas~P. Kennedy.
\newblock The {G}alton-{W}atson process conditioned on the total progeny.
\newblock {\em J. Appl. Probability}, 12(4):800--806, 1975.

\bibitem[Kes86]{Kesten:MR871905}
Harry Kesten.
\newblock Subdiffusive behavior of random walk on a random cluster.
\newblock {\em Ann. Inst. H. Poincar\'e Probab. Statist.}, 22(4):425--487,
  1986.

\bibitem[KS88]{Karatzas--Shreve}
Ioannis Karatzas and Steven~E. Shreve.
\newblock {\em Brownian motion and stochastic calculus}, volume 113 of {\em
  Graduate Texts in Mathematics}.
\newblock Springer-Verlag, New York, 1988.

\bibitem[LG99]{LeGall:MR1714707}
Jean-Fran{\c{c}}ois Le~Gall.
\newblock {\em Spatial branching processes, random snakes and partial
  differential equations}.
\newblock Lectures in Mathematics ETH Z\"urich. Birkh\"auser Verlag, Basel,
  1999.

\bibitem[Pro04]{Protter:MR2020294}
Philip~E. Protter.
\newblock {\em Stochastic integration and differential equations}, volume~21 of
  {\em Applications of Mathematics (New York)}.
\newblock Springer-Verlag, Berlin, second edition, 2004.
\newblock Stochastic Modelling and Applied Probability.

\bibitem[Shi96]{Shiryaev:MR1368405}
A.~N. Shiryaev.
\newblock {\em Probability}, volume~95 of {\em Graduate Texts in Mathematics}.
\newblock Springer-Verlag, New York, second edition, 1996.
\newblock Translated from the first (1980) Russian edition by R. P. Boas.

\bibitem[Wal86]{Walsh-SPDE:MR876085}
John~B. Walsh.
\newblock An introduction to stochastic partial differential equations.
\newblock In {\em \'{E}cole d'\'et\'e de probabilit\'es de {S}aint-{F}lour,
  {XIV}---1984}, volume 1180 of {\em Lecture Notes in Math.}, pages 265--439.
  Springer, Berlin, 1986.

\end{thebibliography}

\end{document}